\documentclass[12pt, reqno]{amsart}

\usepackage[dvipsnames]{xcolor}
\usepackage[colorlinks,citecolor=ForestGreen,urlcolor=blue]{hyperref}
\usepackage{amsmath, amssymb, amsfonts, mathtools}
\usepackage{longtable}
\usepackage{graphicx}
\usepackage{paralist}
\usepackage{csquotes}
\usepackage{dsfont}
\usepackage{pdflscape}
\usepackage{morefloats}
\usepackage{nccfoots}
\usepackage{enumitem, calc} 
\usepackage{tikz}
\usepackage[round]{natbib}
\usepackage{booktabs}
\usepackage{makecell} 
\setcellgapes{2pt}
\usepackage{xcolor}
\allowdisplaybreaks

\usepackage{setspace}
\newcommand{\Var}{{\operatorname{Var}}}
\newcommand{\Cov}{{\operatorname{Cov}}}
\newcommand{\Cor}{{\operatorname{Corr}}}

\newcommand{\bfe}{\mathbf{e}}
\newcommand{\bfp}{\mathbf{p}}
\newcommand{\bfu}{\mathbf{u}}
\newcommand{\bfx}{\mathbf{x}}
\newcommand{\bfX}{\mathbf{X}}
\newcommand{\bfy}{\mathbf{y}}
\newcommand{\bfY}{\mathbf{Y}}
\newcommand{\bbr}{\mathbb{R}}
\newcommand{\bbz}{\mathbb{Z}}
\newcommand{\bbn}{\mathbb{N}}
\newcommand{\bbe}{\mathbb{E}}
\newcommand{\bbp}{\mathbb{P}}
\newcommand{\ca}{\mathcal{A}}
\newcommand{\bfZ}{\mathbf{Z}}
\newcommand{\bfDelta}{\mbox{\boldmath $\Delta$}}
\newcommand{\bfDeltai}{\mbox{\boldmath\scriptsize $\Delta$}}

\newcommand{\abs}[1]{\left| #1 \right|}
\newcommand{\norm}[1]{\left\| #1 \right\|}

\usepackage{dsfont}
\newcommand{\indicator}[1]{\mathds{1}_{\{ #1 \}}}

\newcommand{\brackets}[1]{\left( #1 \right)}

\newcommand{\ignore}[1]{}

\theoremstyle{plain}
\newtheorem{theorem}{Theorem}[section]
\newtheorem{corollary}[theorem]{Corollary}
\newtheorem{lemma}[theorem]{Lemma}
\newtheorem{proposition}[theorem]{Proposition}

\theoremstyle{definition}
\newtheorem{definition}[theorem]{Definition}

\theoremstyle{remark}
\newtheorem{remark}[theorem]{Remark}
\newtheorem{example}[theorem]{Example}

\makeatletter
\newcommand*{\defeq}{\mathrel{\rlap{%
                     \raisebox{0.3ex}{$\m@th\cdot$}}%
                     \raisebox{-0.3ex}{$\m@th\cdot$}}%
                     =}
\makeatother

\makeatletter
\newcommand*{\eqdef}{=\mathrel{\rlap{%
                     \raisebox{0.3ex}{$\m@th\cdot$}}%
                     \raisebox{-0.3ex}{$\m@th\cdot$}}%
                     }
\makeatother
\newcommand\scalemath[2]{\scalebox{#1}{\mbox{\ensuremath{\displaystyle #2}}}}

\usepackage{listings}
\begin{document}
\title[]
{The Joint Asymptotic Distribution of Entropy and Complexity} 
\author[A. Silbernagel]{Angelika Silbernagel}
\address{Department of Mathematics and Statistics, Helmut Schmidt University, 22043 Hamburg, Germany}
\email{\href{mailto:silbernagel@hsu-hh.de}{\nolinkurl{silbernagel@hsu-hh.de}}}

\author[C.H. Wei\ss]{Christian H. Wei\ss}
\address{Department of Mathematics and Statistics, Helmut Schmidt University, 22043 Hamburg, Germany. ORCID: \href{https://orcid.org/0000-0001-8739-6631}{\nolinkurl{0000-0001-8739-6631}}.}
\email{\href{mailto:weissc@hsu-hh.de}{\nolinkurl{weissc@hsu-hh.de}}}
%\urladdr{https://orcid.org/0000-0001-8739-6631}

%\thanks{Research support}
\today


\begin{abstract}
We derive the asymptotic distribution of ordinal-pattern frequencies under weak dependence conditions and investigate the long-run covariance matrix not only analytically for moving-average, Gaussian, and the novel generalized coin-tossing processes, but also approximately by a simulation-based approach. Then, we deduce the asymptotic distribution of the entropy--complexity pair, which emerged as a popular tool for summarizing the time-series dynamics. Here, we make the necessary distinction between a uniform and a non-uniform ordinal pattern distribution and, thus, obtain two different limit theorems. On this basis, we consider a test for serial dependence and check its finite--sample performance. Moreover, we use our asymptotic results to approximate the estimation uncertainty of entropy--complexity pairs.
\end{abstract}
\maketitle

%\vspace*{-7mm} %%%%%%%%%%%%%%%%%%%%%%%%%%%%%%%%%
 \section{Introduction}\label{intro}
About a quarter century ago, in the seminal paper by \citet{bandt02}, a new look at time series $x_1,\ldots,x_T$ was proposed. In order to analyze the dynamic structure (serial dependence) of the data-generating process (DGP) $(X_t)_{t\in\bbz}$, assumed to be real-valued and continuously distributed, it was suggested to transform the given time series into a series of ``permutation types'' (nowadays ``ordinal patterns'', OPs for short), and to conclude from the marginal distribution of the OP-series (a vector of relative frequencies) onto the serial structure of the original process $(X_t)$. As the main advantages of their approach, \citet[p.~1]{bandt02} emphasize its ``simplicity, extremely fast calculation, robustness, and invariance with respect to nonlinear monotonous transformations''. 
More precisely, if the given time series $x_1,\ldots,x_T$ has been transformed into the OPs $\Pi_1,\ldots,\Pi_n$, and if~$\hat\bfp$ is the resulting vector of OP-frequencies, then \citet{bandt02} suggested to compute the Shannon entropy from~$\hat\bfp$, referred to as the ``permutation entropy'', in order to gain insights into the serial properties of $(X_t)$. A few years later, \citet{lamberti04} and \citet{rosso07} suggested to complement the permutation entropy by an additional ``complexity'' measure computed from~$\hat\bfp$, which provides an ``important additional piece of information regarding the peculiarities of a probability distribution, not already carried by the entropy'' \citep[p.~2]{rosso07}. Hence, instead of solely looking at a time series' permutation entropy, the pair of complexity and entropy is used to summarize the essential time-series dynamics. Here, it is important to note that the range of possible entropy--complexity pairs is compact and can be bounded by two curves, see \citet{martin06} for their computation. This motivates to plot the actual entropy--complexity pair within this region, i.e., to visualize the actual time-series dynamics relative to all possible entropy--complexity pairs. The resulting graphical tool, a two-dimensional diagram, is called the ``entropy--complexity plane''. Over the years, it has been successfully applied in various contexts, see \citet{zunino10,ribeiro12,zunino16,sigaki18,chagas22} for some examples. The required curves of minimal and maximal complexity values \citep{martin06} can be computed by using the R~package \href{https://cran.r-project.org/package=statcomp}{\nolinkurl{statcomp}}.

\smallskip
If having estimated the entropy--complexity pair from the given time-series data, and if one wants to infer on the true entropy--complexity pair of the underlying DGP, it is important to consider the actual estimation uncertainty. While it is hardly possible to derive the exact (joint) sampling distribution of entropy and complexity, a promising approach would be to derive an asymptotic expression instead, and to use this expression for approximating the estimation uncertainty. To our knowledge, however, such a joint asymptotic distribution for entropy and complexity has not been derived so far. We are only aware of the approximate solution of \citet{rey25}, where the so-called ``multinomial model'' is used, which ignores the serial dependence within the OP~series. Therefore, the main task of the present paper is the derivation of the true joint asymptotic distribution of entropy and complexity, where we distinguish between the cases of a uniform or non-uniform OP~distribution, and where we consider a wide range of possible DGPs. 

\smallskip
This article is structured as follows: First, we derive the asymptotic distribution of OP frequencies from two perspectives, namely from the perspective of the observed process as well as from the perspective of the resulting OP series, where the latter perspective is relevant if nothing is known about the original process (Section~\ref{section: OP frequencies distribution}). In Section~\ref{section: autocovariance matrix}, we take a closer look at the long-run covariance matrix from the asymptotic results derived before, and we obtain closed-form expressions for important special cases (including the newly proposed generalized coin-tossing process). Moreover, we suggest a simulation-based approach for cases where nothing is known about the actual DGP. In Section~\ref{section: statistics distribution}, we derive the asymptotic distribution of the entropy--complexity pair, where we distinguish between the cases of a uniform and non-uniform OP-distribution. Finally, we discuss some applications of our asymptotic results in Section~\ref{section: applications}. In particular, we consider a test for serial dependence and investigate its finite--sample performance. Furthermore, we suggest a method to quantify the estimation uncertainty of entropy--complexity pairs as well as methods to visualize it. We conclude in Section~\ref{Conclusions} by pointing out some relevant directions for future research. Appendix~\ref{section: approximating functionals and mixing conditions} contains the relevant material on mixing conditions and approximating functionals. All proofs are provided in Appendix~\ref{section: proofs}. In Appendix~\ref{section: MA(q) m=3}, we discuss the long-run covariance matrix with regard to moving-average (MA) processes and OPs of length 3 in more detail.

%%%%%%%%%%%%%%%%%%%%%%%%%%%%%%%%%%%%%%%%%%%%%%%%%%%%
\section{Asymptotic Distribution of the OP Frequencies}
\label{section: OP frequencies distribution}
 
First, let us formally introduce the notion of OPs.

\begin{definition}
    For $m \geq 2$, let $S_m$ denote the set of permutations of $1, \dots, m$, i.e., the tuples which contain each of those numbers exactly once. Let $\bfx = (x_1, \dots, x_m)^\top \in \bbr^m$ be an $m$-dimensional vector. We call the tuple $\pi = (\pi^{(1)}, \dots, \pi^{(m)})$ the \emph{OP of $\bfx$} if the condition
    \begin{equation*}
        \pi^{(j)} < \pi^{(k)} \ \Longleftrightarrow \ x_j < x_k \ \text{ (or } x_j=x_k \text{ for } j < k)
    \end{equation*}
    is satisfied for all $j,k\in \{1,2,\ldots,m\}$.
\end{definition}

There are various ways to deal with ties in the definition of OPs, see \citet{schnurrfischerties} and \citet{schn_sil_23} for a short overview. Here, we assume that the distributions under consideration are continuous, so the probability of ties is zero.

\begin{table}[t]
    \label{tab: lexicographic order}
    \caption{Lexicographic ordering of OPs of length $m=3$.}
    \centering
    \begin{tabular}{c|cccccc}
        \toprule
        $i$ & 1 & 2 & 3 & 4 & 5 & 6 \\
        \midrule
        $\pi_i$ & (1,2,3) & (1, 3, 2) & (2, 1, 3) & (2, 3, 1) & (3, 1, 2) & (3, 2, 1) \\
        \bottomrule
    \end{tabular}
\end{table}

\smallskip
There exist $m!$ OPs of length $m$, which we denote by $\pi_1, \dots, \pi_{m!}$ sorted according to the lexicographic order. For $m=2$, this is $\pi_1 = (1,2)$ and $\pi_2=(2,1)$. The numbering for $m=3$ is presented in Table~\ref{tab: lexicographic order}. 
The OP of $\bfx \in \bbr^m$ is also uniquely determined by the vector of increments $\bfDelta_{\bfx} = (x_2-x_1, \dots, x_m-x_{m-1})^\top\in\bbr^{m-1}$. 

\medskip
The minimum assumption necessary to derive reasonable limit theorems for the OP frequency distribution is \emph{order stationarity}, that is, the probability for each OP does not change over time (see, e.g., \citealp{bandt25}). From a time-series point of view, however, it might be more intuitive to impose a stationarity assumption on the observed process $(X_t)_{t \in \bbz}$ instead of the series of observed OPs. A sufficient condition here is that of stationary increments $\Delta_t = X_t - X_{t-1}$ instead of stationarity of the observed process itself. Note that the stationarity of a process implies stationary increments which, in turn, implies order stationarity.

\smallskip
In the present paper, we consider the asymptotic distribution of OP frequencies from two perspectives, namely from the perspective of the observed process $(X_t)_{t \in \bbz}$ as well as from the perspective of the resulting OP-series $(\Pi_t)_{t \in \bbz}$. On this basis, we formulate suitable conditions that lead to the desired limit theorems.

\subsection{Time-Series Perspective}
\label{subsection: time series approach}

Let $(X_t)_{t\in\bbz}$ denote the observed real-valued process and suppose it is increment stationary, i.e., the process $(\Delta_t)_{t\in\bbz}$ with $\Delta_t = X_t - X_{t-1}$ is stationary.
As we are interested in the OPs stemming from $(X_t)_{t\in\bbz}$, it is convenient to consider the process $(\bfDelta_t)_{t\in\bbz}$ defined by
\[
    \bfDelta_t:={(\Delta_{t+1},\ldots, \Delta_{t+m-1})^\top} = (X_{t+1} - X_t, \dots, X_{t+m-1} - X_{t+m-2})^\top,\quad t \in \bbz,
\]
which is obtained from $(\Delta_t)_{t\in\bbz}$ by a moving-window approach. 
Note that if $(\Delta_t)_{t\in \bbz}$ is stationary, then $(\bfDelta_t)_{t\in \bbz}$ is stationary as well.

\smallskip
Let $\Pi : \bbr^{m} \to S_m$ denote the function which assigns each $m$-dimensional vector its OP, while $\Pi_{\bfDeltai} : \bbr^{m-1} \to S_m$ denotes the function which assigns each $(m-1)$-dimensional \emph{vector of increments} its OP. Writing $\bfX_t = (X_t, \dots, X_{t+m-1})^\top$ for each $t\in\bbz$, let us point out that $\Pi_t := \Pi(\bfX_t) = \Pi_{\bfDeltai}(\bfDelta_t)$ holds. 
Furthermore, let $\bfp = (p_1, \dots, p_{m!})^\top$ denote the vector consisting of the OP probabilities, i.e., $\bfp$ denotes the probability mass function (pmf) vector of the OPs. A natural estimator for each probability $p_i$ is given by the relative frequency $\hat{p}_i = \tfrac{1}{n} \sum^n_{t=1} \indicator{\Pi(\bfX_t)=\pi_i} = \tfrac{1}{n} \sum^n_{t=1} \indicator{\Pi_{\bfDeltai}(\bfDeltai_t)=\pi_i}$. Hence, defining a kernel $h : \bbr^{m-1} \to \bbr^{m!}$ by
\begin{equation}\label{eq: kernel}
    h(\bfy):=(\indicator{\Pi_{\bfDeltai}(\bfy)=\pi_1}, \dots, \indicator{\Pi_{\bfDeltai}(\bfy)=\pi_{m!}})^\top,
\end{equation}
we can estimate $\bfp$ by $\hat{\bfp} = (\hat{p}_1, \dots, \hat{p}_{m!})^\top = \tfrac{1}{n} \sum^n_{t=1} h(\bfDelta_t)$.

\smallskip
In order to derive the asymptotic distribution of $\hat{\bfp}$ based on the process $(X_t)_{t\in\bbz}$, we shall make use of the concept of ``r-approximating functionals'' which is sometimes referred to by ``near-epoch dependence'' in the literature. The idea is that even though a process does not satisfy a certain mixing condition itself, it may be approximated almost entirely by the ``near epoch'' of a ``background process'' $(Z_t)_{t\in\bbz}$ which satisfies the required mixing condition. Here, we resort to absolute regularity as the relevant mixing condition. A short overview providing all necessary definitions is given in Apppendix~\ref{section: approximating functionals and mixing conditions}.

\smallskip
\citet{schn_sil_mar_24} have shown that if the process $(\Delta_t)_{t\in \bbz}$ is an $r$-approximating functional with approximating constants $(a_k)_{k \in \bbn_0}$ of size $-\lambda$, then $(\bfDelta_t)_{t\in\bbz}$ also satisfies the $r$-approximating condition for $k \geq m-1$ with approximating constants $(a^\prime_k)_{k \in \bbn_0}$ defined by $a^\prime_k = 2 \sum^{(m-1)-1}_{i=0} a_{k-i}$. Note that this can be improved to $k \geq (m-1)-1=m-2$, since there is no argument in the authors' proof that contradicts such a generalization.
In order to derive the asymptotic distribution of $\hat{\bfp}$, first we show that our kernel satisfies the technical condition of ``$p$-continuity'' (see Definition~\ref{definition: p-continuity}). The proof is provided in Appendix~\ref{proof prop: p-continuous}.

\begin{proposition} \label{prop:p-continuous}
Let $p\geq 1$. Then, the kernel $h$ as defined in \eqref{eq: kernel} is $p$-continuous with respect to the distribution of $\bfDelta_0$.
\end{proposition}

The limit distribution of the pmf vector is now a consequence of \citet[Theorem~B.5]{schn_sil_mar_24}, which is a central limit theorem (CLT) for partial sums of functionals of absolutely regular processes. For the readers' convenience, we have stated it in the multivariate version in Theorem~\ref{theorem: Theorem 4}.

\begin{theorem}\label{theorem: asymp distribution OP probs case 1}
    Let $(\Delta_t)_{t\in\bbn}$ be a 1-approximating functional with approximating constants $(a_k)_{k \in \bbn_0}$ of an absolutely regular process with mixing coefficients $(\beta_k)_{k\in\bbn_0}$ such that 
    \begin{equation*}
        \sum^{\infty}_{k=0} k^2 \brackets{\sqrt{a_k} + \beta_k} < \infty.
    \end{equation*}
    Assume moreover that
    \begin{equation}
    \label{eq: asymp distribution OP probs case 1 phi}
        \sum^{\infty}_{k=m-1} k^2 \phi\left(2 \sqrt{\sum^{m-1}_{i=0} a_{k-i}}\right) < \infty,
    \end{equation}
    where $\phi$ is defined in \eqref{eq: phi specified}.
    Then, as $n \to \infty$,
    \begin{equation*}
        \sqrt{n} (\hat{\bfp}-\bfp) \xrightarrow{d} N(\mathbf{0}, \Sigma),
    \end{equation*}
    where the long-run covariance matrix $\Sigma = (\sigma_{ij})_{1 \leq i,j \leq m!}$, defined by
    \begin{equation}\label{eq: long-run covariance matrix}
        % \sigma_{ij} = p_i (\delta_{ij} - p_j) + \sum^\infty_{k=1} \Bigl(\bbp(\Pi(\bfX_0) = \pi_i, \Pi(\bfX_k)=\pi_j) + \bbp(\Pi(\bfX_0) = \pi_j, \Pi(\bfX_k)=\pi_i) - 2 p_i p_j\Bigr)
        \sigma_{ij} = p_i (\delta_{ij} - p_j) + \sum^\infty_{k=1} \Bigl(p_{ij}(k) + p_{ji}(k) - 2 p_i p_j\Bigr)
    \end{equation}
    with $p_{ij}(k) = \bbp(\Pi(\bfX_0) = \pi_i, \Pi(\bfX_k)=\pi_j)$, 
    converges absolutely with $\delta_{ij} = \indicator{i=j}$ denoting the Kronecker delta.
\end{theorem}

The proof is provided in Appendix~\ref{proof theorem: asymp distribution OP probs case 1}.
Note that the CLT derived by \citet[Section~2]{wei_22} is a special case of Theorem~\ref{theorem: asymp distribution OP probs case 1}.
Also compare with Theorem~3 (and especially Lemma~1~(i)) from \citet{schn_deh_17}, where the additional assumption of $X_i-X_1$ being Lipschitz-continuous for all $i \in \{1, \dots, m\}$ is imposed.

\smallskip
The conditions of Theorem~\ref{theorem: asymp distribution OP probs case 1} could have been posed directly onto the observed process $(X_t)_{t\in\bbz}$ instead of the increment process $(\Delta_t)_{t\in\bbz}$. However, proceeding like in Theorem~\ref{theorem: asymp distribution OP probs case 1}, we obtain a more general result, see Appendix~\ref{appendix remark}.

\begin{remark}
    \label{remark: asymp distribution OP probs case 1}
    \begin{enumerate}[label=(\alph*)]
        \item The assumptions of Theorem~\ref{theorem: asymp distribution OP probs case 1} are automatically satisfied if the respective processes are $q$-dependent, because these are absolutely regular with coefficients $\beta_k=0$ for $k>q$ due to independence of the sub-$\sigma$-fields under consideration. As the process is already absolutely regular, no approximation is needed, so $a_k \equiv 0$. Hence, by the definition of $\phi$ in \eqref{eq: phi specified}, \eqref{eq: asymp distribution OP probs case 1 phi} reduces to 
        \begin{equation}
            \label{eq: phi series example}
            \sum^\infty_{k=m-1} k^2 \phi(0) = \sum^\infty_{k=m-1} k^2 m! \cdot \bbp(\norm{\bfY-\bfY^\prime}=0) = 0,
        \end{equation}
        where $\bfY$ and $\bfY^\prime$ denote two $(m-1)$-dimensional vectors with the same distribution as $\bfDelta_0$. Recall that $\bfDelta_0$ is continuously distributed. Therefore, all assumptions of Theorem~\ref{theorem: asymp distribution OP probs case 1} are met. 
        \item Similar arguments show that autoregressive moving-average (ARMA) processes defined on innovations, the probability law of which is absolutely continuous with regard to the Lebesgue measure on $\bbr^d$, fulfill the conditions of Theorem~\ref{theorem: asymp distribution OP probs case 1}: \cite{mokkadem88} has shown that such processes are absolutely regular and that there exists a constant $\rho \in (0,1)$ such that the mixing coefficients $(\beta_k)_{k \geq 0}$ satisfy $\beta_k = \mathcal{O}(\rho^k)$. Note that this convergence rate ensures the required convergence of $\sum^\infty_{k=1} k^2 \beta_k$. Moreover, no approximation to an absolutely regular process is needed, so $a_k \equiv 0$, and hence, \eqref{eq: phi series example} holds true. 
    \end{enumerate}
\end{remark}

%%%%%%%%%%%%%
\subsection{OP-Series Perspective}
\label{subsection: OP-series}
Now, let us turn to the OP-series perspective as suggested by \citet{bandt25}. In this case, we solely focus on the resulting OP~series $(\Pi_t)_{t\in\bbz}$ and its stochastic properties, but the underlying original DGP $(X_t)_{t\in\bbz}$ is not further specified. It suffices to assume that $(X_t)_{t\in\bbz}$ is an at least ordinal process consisting of certain ``objects'' $X_t$, which are not necessarily numbers but still allow for the required strict ordering (i.e., in order to avoid dealing with ties, we assume that the probability of ties is zero). An example for such a process is the coin-tossing model introduced by \citet{bandt25}, which is discussed in more detail in Section~\ref{subsection: coin-tossing order}. 
Let $(\Pi_t)_{t\in\bbz}$ denote the process consisting of the observed OPs $\Pi_t := \Pi(\bfX_t)$ at time $t$. Instead of posing a stationarity assumption to $(X_t)_{t\in\bbz}$, we assume that $(X_t)_{t\in\bbz}$ is ``order stationary'', that is, we require the stationarity of $(\Pi_t)_{t\in\bbz}$.
Moreover, suppose that $(\Pi_t)_{t\in\bbz}$ constitutes a strongly mixing process.
Here, we choose this weak dependence property in order to not have to deal with the issue of a proper definition for distances between nominal categories (here given by the observed OPs) for instance, which would arise when considering approximating functionals.

\begin{theorem}
    \label{theorem: asymp distribution OP probs case 2}
    Let $(\Pi_t)_{t\in\bbn}$ be stationary and strongly mixing with mixing coefficients $(\alpha_k)_{k\in\bbn_0}$ such that 
    \begin{equation*}
        \sum^{\infty}_{k=0} \alpha_k < \infty.
    \end{equation*}
    Then, as $n \to \infty$,
    \begin{equation*}
        \sqrt{n} (\hat{\bfp}-\bfp) \xrightarrow{d} N(\mathbf{0}, \Sigma)
    \end{equation*}
    where $\Sigma$ is given in \eqref{eq: long-run covariance matrix}.
\end{theorem}

Defining the function $g : S_m \to \{0,1\}^{m!}$ by
$g(\Pi_t) := (\indicator{\Pi_t = \pi_1}, \dots, \indicator{\Pi_t = \pi_{m!}})^\top$,
we can binarize the observed OP $\Pi_t$ by using $g$. 
Writing $\hat{\bfp} = \tfrac{1}{n} \sum^n_{t=1} g(\Pi_t)$, the claim of Theorem~\ref{theorem: asymp distribution OP probs case 2} follows immediately from \citet[Theorem~18.6.3]{ibragimovlinnik71}. 

\smallskip
Note that the CLTs derived in \citet[Section~4]{sousa22} and \citet[Theorem~1]{elsinger10} are special cases of Theorem~\ref{theorem: asymp distribution OP probs case 2}, since $q$-dependent processes are strongly mixing with coefficients $\alpha_k = 0$ for all $k > q$. While the random walk itself is not stationary, it can be easily shown that the OP process derived from it is not only stationary, but also $(m-2)$-dependent (see \citealp{sousa22}). Hence, it satisfies the assumptions of Theorem~\ref{theorem: asymp distribution OP probs case 2}. 
Another example for a process, where the OP series is $(m-2)$-dependent and therefore satisfies the above assumptions, is the coin-tossing order (see Section~\ref{subsection: coin-tossing order}). 

\smallskip
Moreover, one can show that if the original process $(X_t)_{t\in\bbz}$ is strongly mixing, then the same holds true for the OP process (with coefficients of same size) as $\Pi_t = \Pi(X_t, \dots, X_{t+m-1})$ is just a measurable function evaluated in a finite number of arguments (see \citealp[Theorem~14.1]{davidson}). Some more specific examples are the following.

\begin{example}
    \label{example: asymp distribution OP probs case 2}
    \begin{enumerate}[label=(\alph*)]
        \item Let us consider the ARMA processes from Remark~\ref{remark: asymp distribution OP probs case 1}~(b) again, but now from this new perspective. As such processes are absolutely regular, they also satisfy the strong mixing condition with coefficients $a_k \leq \beta_k = \mathcal{O}(\rho^k)$ where $0<\rho<1$. The convergence of $\sum^\infty_{k=0} \alpha_k$ now follows by direct comparison test. Consequently, the above considerations imply that such ARMA processes also fulfill the conditions of Theorem~\ref{theorem: asymp distribution OP probs case 2}.
        \item The first-order transposed exponential autoregressive (TEAR$(1)$) process $(X_t)_{t\in\bbz}$ is defined by the recursion $X_t=B_t\,X_{t-1}+(1-p)\,\varepsilon_t$ with i.i.d.\ $\varepsilon_t\sim \text{Exp}(1)$ and Bernoulli $B_t$ with $P(B_t=1)=p$. Since $(X_t)_{t\in\bbz}$ is a stationary Markov chain (see \citealp{lawrance81}), the strong mixing coefficients reduce to 
        \[
            \alpha_k = \alpha(\sigma(X_0), \sigma(X_k)) = \sup_{A \in \sigma(X_0), B \in \sigma(X_k)} |\bbp(A \cap B) - \bbp(A) \bbp(B)| \leq p^k 
        \]
        for $k \geq 0$ (see \citealp{bradley05}). 
        The inequality follows by induction as $X_0$ and $X_1$ are independent with probability $1-p$. Hence, $(X_t)_{t\in\bbz}$ is strongly mixing with summable mixing coefficients and, according to the above remark, the OP~process resulting from $(X_t)_{t\in\bbz}$ has the same properties.
    \end{enumerate}
\end{example}

%%%%%%%%%%%%%%%%
\section{Long-run Covariance Matrix of OP~Frequencies}
\label{section: autocovariance matrix}

So far, the existing literature on explicit derivations of the covariance matrix in \eqref{eq: long-run covariance matrix} is very limited. In the practically most relevant case of OP-length $m=3$ \citep[see][]{bandt19}, results are available for a random walk with symmetric noise (``sRW'', see \citealp{sousa22}) and i.i.d.\ sequences (derived independently by \citealp{wei_22} and \citealp{sousa22}), namely
\begin{equation}
\label{Sigma_sRW_iid}
\textstyle
\Sigma_{\textup{sRW}}\ =\ 
\frac{1}{192}\,
\begin{psmallmatrix*}[r]
 60 & -6 & -6 & -6 & -6 & -36 \\
 -6 & 15 & 7 & -9 & -1 & -6 \\
 -6 & 7 & 15 & -1 & -9 & -6 \\
 -6 & -9 & -1 & 15 & 7 & -6 \\
 -6 & -1 & -9 & 7 & 15 & -6 \\
-36 & -6 & -6 & -6 & -6 & 60 \\
\end{psmallmatrix*}
\text{ and}\quad
\Sigma_{\textup{iid}}\ =\ 
\frac{1}{360}\,
\begin{psmallmatrix*}[r]
 46 & -23 & -23  &  7  &  7 & -14 \\
-23 &  28 &  10 & -20 &  -2  &  7 \\
-23 &  10 &  28 &  -2 & -20  &  7 \\
  7 & -20 &  -2 &  28 &  10 & -23 \\
  7 &  -2 & -20 &  10 &  28 & -23 \\
-14  &  7  &  7 & -23 & -23 &  46 \\
\end{psmallmatrix*},
\end{equation}
respectively. Note that the covariance matrices in \eqref{Sigma_sRW_iid} refer to the lexicographic order of OPs according to Table~\ref{tab: lexicographic order}, whereas \citet{wei_22,sousa22} used different order schemes. Moreover, \citet{sousa22} have estimated the covariance matrix for an MA process of order $q$ (``MA($q$)'' for short) for $m=3$ with non-necessarily equal weights. To the best of our knowledge, no further work has been done on the covariance matrix in \eqref{eq: long-run covariance matrix} so far. In what follows, we further reduce this gap by contributing three novel cases, namely an MA($q$) process with equal weights and Gaussian processes (both for OPs of length $m=2$), as well as a generalization of the coin-tossing process introduced by \citet{bandt25}.

\subsection{MA(q) Process}
\label{subsection: MA(q)}

Let $(\varepsilon_t)_{t\in\bbz}$ denote an i.i.d.\ series.
We consider the MA($q$) process with equal weights, that is,
\begin{equation*}
    X_t = \sum^q_{l=0} \varepsilon_{t-l}.
\end{equation*}
For a fixed OP length $m$, this results in an $(m-1+q)$-dependent OP process (see \citealp{sousa22}). Therefore, the sum in \eqref{eq: long-run covariance matrix} breaks off at this point. We provide the covariance matrix for $m=2$. Analogous computations are conceivable for the case $m=3$, but they are more involved and the solution depends on the distribution of the innovations, see Appendix~\ref{section: MA(q) m=3} for more details. The proof of Proposition~\ref{prop: cov MA(q) m2} is provided in Appendix~\ref{proof prop: cov MA(q) m2}.

\begin{proposition}
    \label{prop: cov MA(q) m2}
    For an MA($q$) process with equal weights and OPs of length $m=2$, the covariance matrix from \eqref{eq: long-run covariance matrix} is given by
\begin{equation*}
    \Sigma = \frac{1}{12} \left(\begin{array}{rr}
        1 & -1 \\
        -1 & 1
    \end{array}\right).
\end{equation*}
\end{proposition}

In particular, Proposition~\ref{prop: cov MA(q) m2} states that $\Sigma$ is independent of the order $q$ of the moving-average process as well as from the distribution of the background process. Comparing this matrix to the matrix derived by \citet{wei_22}, we also recognize that for OPs of length $m=2$, the covariance matrices for an i.i.d.\ and an MA($q$) process with equal weights coincide.

\subsection{q-Dependent Gaussian Processes}
\label{subsection: q-dep Gaussian}

For $q\geq1$, let $(X_t)_{t\in\bbz}$ be a $q$-dependent Gaussian process with stationary increments. Then, the increment process $(\Delta_t)_{t\in\bbz}$ is a stationary Gaussian process with zero mean and autocorrelations 
\begin{align*}
    \rho(k):= \Cor(\Delta_0, \Delta_k) &= \frac{\Cov(X_0, X_k) - \Cov(X_{-1}, X_k) - \Cov(X_0, X_{k-1}) + \Cov(X_{-1}, X_{k-1})}{\Var(X_0) - 2\Cov(X_{-1}, X_0) + \Var(X_{-1})}.
\end{align*}
Note that $(\Delta_t)_{t\in\bbz}$ is $(q+1)$-dependent. We consider OPs of length $m=2$, i.e., we have to consider probabilities of the form 
\[
    \bbp(X_0 \delta X_1, X_k \delta^\prime X_{k+1}) = \bbp(0\ \delta\ \Delta_1,\ 0\ \delta^\prime \Delta_{k+1})
\]
for $\delta, \delta^\prime \in \{<, >\}$. It is a well-known fact that 
\[
    \bbp(0 > \Delta_1, 0 > \Delta_{k+1}) = \bbp(0 < \Delta_1, 0 < \Delta_{k+1}) = \frac{1}{4} + \frac{1}{2\pi} \arcsin{\rho(k)},
\]
see \citet[Section~5]{bandt07}. Hence,
\[
    \bbp(0 > \Delta_1, 0 < \Delta_{k+1}) = \bbp(0 < \Delta_1, 0 > \Delta_{k+1}) = \frac{1}{4} - \frac{1}{2\pi} \arcsin{\rho(k)}
\]
such that
\begin{align}
    \Sigma &= \frac{1}{4} \cdot \begin{pmatrix}
        1 & -1 \\ -1 & 1
    \end{pmatrix} + \sum^{q+1}_{k=1} \brackets{2\cdot\brackets{\frac{1}{4} + \begin{pmatrix}
        1 & -1 \\ -1 & 1
    \end{pmatrix} \cdot \frac{1}{2\pi} \arcsin{\rho(k)}}- 2 \cdot \frac{1}{2} \cdot \frac{1}{2}} \nonumber \\
    &= \brackets{\frac{1}{4} + \frac{1}{2\pi}\sum^{q+1}_{k=1} \arcsin{\rho(k)}}\cdot\begin{pmatrix}
        1 & -1 \\ -1 & 1
    \end{pmatrix}. \label{eq: Gaussian process autocovariance matrix}
\end{align}

An example for such a $q$-dependent process is a Gaussian MA$(q)$-process. In Section~\ref{subsection: MA(q)}, we have studied an MA($q$) process with equal weights. Now, we extend our considerations to the case
\begin{equation*}
    X_t = \varepsilon_t + \sum^q_{l=1} \theta_i \varepsilon_{t-l},
\end{equation*}
where $\varepsilon_t \sim N(\mu, \sigma^2)$ is an i.i.d.\ series and $\theta_1, \dots, \theta_q \in \bbr$. Note that $(X_t)_{t\in\bbz}$ is necessarily a stationary Gaussian process with autocovariances
\[
    \Cov(X_0, X_k) = \begin{cases}
        \sigma^2 (\sum^q_{l=k} \theta_l \theta_{l-k}), &0 \leq k \leq q, \\
        0, &k > q,
    \end{cases}
\]
where $\theta_0:=1$. Accordingly,
\begin{align*}
    \rho(k) &= \frac{2 \Cov(X_0, X_k) - \Cov(X_0, X_{k+1}) - \Cov(X_0, X_{k-1})}{2(\Var(X_0) - \Cov(X_0, X_1))} \\
    &=\frac{2 \Cov(X_0, X_k) - \Cov(X_0, X_{k+1}) - \Cov(X_0, X_{k-1})}{2\sigma^2(\sum^q_{l=0} \theta_l^2 - \sum^q_{l=1} \theta_l \theta_{l-1})}
\end{align*}
such that the autocovariance matrix in \eqref{eq: long-run covariance matrix} can be computed explicitly using \eqref{eq: Gaussian process autocovariance matrix}.
However, let us emphasize that this is restricted by the assumption of $(\varepsilon_t)_{t\in\bbz}$ being Gaussian and OPs of length $m=2$. There is no hope of obtaining explicit formulas for $m\geq3$ as this would require a formula for $\bbp(Y_1>0, Y_2>0, Y_3>0, Y_4>0)$, where $(Y_1, \dots, Y_4)$ has a zero mean Gaussian distribution; also see \citet[Section~5]{bandt07}.

\subsection{Coin-Tossing Order}
\label{subsection: coin-tossing order}

Lastly, let us discuss the coin-tossing order as introduced by \cite{bandt25}, 
where a random order is defined via repeated coin tossing in the following way. For two ``objects'' $X_t$ and $X_{t+1}$ of the underlying original DGP, we throw a 
coin to decide whether $X_t < X_{t+1}$ or $X_t > X_{t+1}$, and we do the same for determining the order between $X_{t+1}$ and $X_{t+2}$. These coin tosses are independent of each other. The interesting part is now the order between $X_t$ and $X_{t+2}$. If it is already fixed by transitivity, which would be the case for $X_t < X_{t+1} < X_{t+2}$ or $X_t > X_{t+1} > X_{t+2}$, then we turn to $X_{t+3}$. Otherwise, we flip another (independent) 
coin to determine the ordering between $X_t$ and $X_{t+2}$ first, before moving forward to $X_{t+3}$. The next step would be to compare $X_{t+3}$ with all objects before unless the order is already fixed by previous comparisons and transitivity, and so on. Note that coin-tossing is always needed for deciding the order of $X_t$ and $X_{t-1}$. Furthermore, the coin-tossing for $X_t$ and $X_{t-j}$, $j \in \bbn$, is carried out regardless of the order of $\ldots, X_{t-j-1}, X_{t-j}$, so it is independent of all coin tosses with regard to pairs of objects in $\{X_s, {s\leq t-j}\}$. Considering OPs of length $m$, the OP process derived from the coin-tossing order is thus $(m-2)$-dependent.

\smallskip
For $i<j$, let $C_{ij}^+$ and $C_{ij}^-$ denote the events $X_i < X_j$ and $X_i > X_j$, respectively. To shorten notation, we leave out the subscripts whenever we want to indicate increasing or decreasing ordering between two arbitrary (but distinct) elements $X_i$ and $X_j$. \cite{bandt25} has introduced the coin-tossing model for fair coins, i.e., where the probabilities for increasing and decreasing orderings are given by $\bbp(C^+)=\bbp(C^-)=\tfrac{1}{2}$. We extend this model to coin tosses where $\bbp(C^+) = p \in (0,1)$ and thus, $\bbp(C^-) = 1-p =:q$ arbitrarily, leading to the ``generalized coin-tossing (GCT) process''. 
First, we consider the case of OPs of length $m=2$. Each OP is determined by one coin toss, so the probability for the increasing and decreasing OP is given by $p$ and $1-p$, respectively. Then, the $(m-2)$-dependence of the OP process already yields
\[
    \Sigma = p(1-p) \cdot\begin{pmatrix}
        1 & -1 \\ -1 & 1
    \end{pmatrix}
\]
as the covariance matrix from \eqref{eq: long-run covariance matrix}.
Next, we turn to the practically more relevant case $m=3$.

\begin{proposition}
\label{prop: Sigma_ct}
    Consider the GCT process with $\bbp(C^+) = p \in (0,1)$ and $q=1-p$. For OPs of length $m=3$, 
    the marginal distribution equals
    \begin{equation}
    \label{marg_ct}
\bfp = \big(p^2,\ p^2 q,\ p^2 q,\ p q^2,\ p q^2,\ q^2\big)^\top,
    \end{equation}
    and the covariance matrix from \eqref{eq: long-run covariance matrix} is given by 
    \begin{equation}
    \label{Sigma_ct}
    \Sigma = 
    \scalemath{0.83}{\begin{pmatrix}
        p^2(1+2p-3p^2) &p^3q(1-3p) &p^3q(1-3p) &p^2q^2(1-3p) &p^2q^2(1-3p) &-3p^2q^2 \\
        p^3q(1-3p) &p^2q(1-3p^2q) &p^3q(1-3pq) &-3p^3q^3 &p^2q^2(1-3pq) &p^2q^2(1-3q) \\
        p^3q(1-3p) &p^3q(1-3pq) &p^2q(1-3p^2q) &p^2q^2(1-3pq) &-3p^3q^3 &p^2q^2(1-3q) \\
        p^2q^2(1-3p) &-3p^3q^3 &p^2q^2(1-3pq) &pq^2(1-3pq^2) &pq^3(1-3pq) &pq^3(3p-2) \\
        p^2q^2(1-3p) &p^2q^2(1-3pq) &-3p^3q^3 &pq^3(1-3pq) &pq^2(1-3pq^2) &pq^3(3p-2) \\
        -3p^2q^2 &p^2q^2(1-3q) &p^2q^2(1-3q) &pq^3(3p-2) &pq^3(3p-2) &pq^2(4-3p)
    \end{pmatrix}}.
    \end{equation}
\end{proposition}

The proof is provided in Appendix~\ref{proof prop: Sigma_ct}.
It is possible to determine the covariance matrices for $m>3$ analogously, though transition probabilities have to be determined for $1 \leq k \leq m-2$ then. 
As a direct consequence of Proposition~\ref{prop: Sigma_ct}, we obtain the following corollary.

\begin{corollary}
\label{prop: Sigma_ct fair coin}
    Consider the GCT process with $\bbp(C^+) = \bbp(C^-) = \tfrac{1}{2}$ (fair coins). For OPs of length $m=3$, 
    the marginal distribution equals $\bfp = \big(\frac{1}{4}, \frac{1}{8}, \frac{1}{8}, \frac{1}{8}, \frac{1}{8}, \frac{1}{4}\big)^\top$, 
    and the covariance matrix from \eqref{eq: long-run covariance matrix} is given by
    \begin{equation}
    \label{Sigma_ct fair coin}
    \Sigma = \frac{1}{64} \scalemath{0.85}{\begin{pmatrix}
        20 & -2 & -2 & -2 & -2 & -12 \\
        -2 & 5 & 1 & -3 & 1 & -2 \\
        -2 & 1 & 5 & 1  & -3 & -2 \\
        -2 & -3 & 1 & 5 & 1 & -2 \\
        -2 & 1 & -3 & 1 & 5 & -2 \\
        -12 & -2 & -2 & -2 & -2 & 20
    \end{pmatrix}}.
    \end{equation}
\end{corollary}

It is worth pointing out that the covariance matrix \eqref{Sigma_ct fair coin} for the fair coin-tossing process with $m=3$ and the one for the symmetric random walk in \eqref{Sigma_sRW_iid} differ from each other only at few places (in the upper triangle, these are the entries $\sigma_{23}$, $\sigma_{25}$, $\sigma_{34}$, and $\sigma_{45}$). The close relationship between these two kinds of process was already pointed out by \citet{bandt25}.

\subsection{Approximate Computation of Covariance Matrix}
\label{Approximate Computation of Covariance Matrix}
As outlined before, there are only a few scenarios where the exact expression of the covariance matrix \eqref{eq: long-run covariance matrix} is known. In order to obtain an at least approximate expression in other cases, we propose a simulation-based approach together with classical methods for covariance estimation. More precisely, one first simulates a ``very long'' time series from either the original DGP $(X_t)_{t\in\bbz}$ and computes the resulting OP-series afterwards (``time-series perspective'' from Section~\ref{subsection: time series approach}), or directly from the OP-process $(\Pi_t)_{t\in\bbz}$ (``OP-series perspective'' from Section~\ref{subsection: OP-series}, e.g., like in Section~\ref{subsection: coin-tossing order}). Here, we considered the length~$n=10^8$ for the final OP-series as ``very long''. Then, in the second step, we consider the binary vectors $\bfZ_t$ derived from~$\Pi_t$ via $Z_{t,i}=\indicator{\Pi_t=\pi_i}$ for $i=1,\ldots,m!$, recall \eqref{eq: kernel}, and apply a procedure for estimating the covariance matrix \eqref{eq: long-run covariance matrix} from $\bfZ_1, \ldots, \bfZ_n$. Here, we used the classical approach of \citet{newey87},
\begin{equation}
\label{NWestimation}
\textstyle
\widehat{\Sigma}\ =\ \widehat{\Omega}_0 + \sum_{k=1}^u w(k,u)\,\big(\widehat{\Omega}_k + \widehat{\Omega}_k^\top\big)
\end{equation}
together with different rules for the truncation parameter~$u$ as well as different weighting schemes~$w(k,u)$, see \citet{lazarus18} for an overview. In \eqref{NWestimation}, $\widehat{\Omega}_k$ abbreviates the sample covariance matrix for the lag-$k$ pairs $(\bfZ_t, \bfZ_{t+k})$. In our experiments, we observed that the Bartlett weights $w(k,u)=1-k/u$ together with the ``textbook rule'' $u=\lfloor0.75\,n^{1/3}\rfloor$ \citep[p.~542]{lazarus18}, or the unit weights $w(k,u)=1$ together with $u=\lfloor n^{1/5}\rfloor$ or  $u=\lfloor 0.75\,n^{1/3}\rfloor$, lead to reasonable results (for $n=10^8$, we get $u=39$ and $u=348$, respectively). Here, the success of the unit weights (even with rather low~$u$) is explained by the fact that often, the OP-series has a very short dependence, i.e., only a few summands in \eqref{eq: long-run covariance matrix} contribute notably to the true covariance matrix~$\Sigma$. 

\smallskip
In order to illustrate the approximation quality of the methods ``Bartlett with $u=348$'' (B348), ``unit with $u=39$'' (U39), and ``unit with $u=348$'' (U348), we compare the obtained approximate values with the true covariance matrix in the aforementioned three cases for $m=3$, namely for an i.i.d.\ DGP like in \eqref{Sigma_sRW_iid},
$$
\begin{array}{l@{\ }ll@{\ }l}
\text{\footnotesize true:} & \begin{psmallmatrix*}[r]
0.128 & -0.064 & -0.064 & 0.019 & 0.019 & -0.039 \\
-0.064 & 0.078 & 0.028 & -0.056 & -0.006 & 0.019 \\
-0.064 & 0.028 & 0.078 & -0.006 & -0.056 & 0.019 \\
0.019 & -0.056 & -0.006 & 0.078 & 0.028 & -0.064 \\
0.019 & -0.006 & -0.056 & 0.028 & 0.078 & -0.064 \\
-0.039 & 0.019 & 0.019 & -0.064 & -0.064 & 0.128 \\
\end{psmallmatrix*},
& \text{\footnotesize B348:} & 
\begin{psmallmatrix*}[r]
0.128 & -0.063 & -0.064 & 0.019 & 0.019 & -0.039 \\
-0.063 & 0.078 & 0.028 & -0.055 & -0.006 & 0.019 \\
-0.064 & 0.028 & 0.078 & -0.006 & -0.056 & 0.019 \\
0.019 & -0.055 & -0.006 & 0.078 & 0.028 & -0.063 \\
0.019 & -0.006 & -0.056 & 0.028 & 0.078 & -0.064 \\
-0.039 & 0.019 & 0.019 & -0.063 & -0.064 & 0.128 \\
\end{psmallmatrix*},
\\[-1ex]
\\
\text{\footnotesize U39:} & 
\begin{psmallmatrix*}[r]
0.128 & -0.064 & -0.064 & 0.019 & 0.020 & -0.039 \\
-0.064 & 0.078 & 0.028 & -0.055 & -0.006 & 0.019 \\
-0.064 & 0.028 & 0.078 & -0.006 & -0.056 & 0.020 \\
0.019 & -0.055 & -0.006 & 0.078 & 0.028 & -0.064 \\
0.020 & -0.006 & -0.056 & 0.028 & 0.078 & -0.064 \\
-0.039 & 0.019 & 0.020 & -0.064 & -0.064 & 0.128 \\
\end{psmallmatrix*},
& \text{\footnotesize U348:} & 
\begin{psmallmatrix*}[r]
0.127 & -0.064 & -0.063 & 0.020 & 0.019 & -0.039 \\
-0.064 & 0.078 & 0.027 & -0.056 & -0.005 & 0.020 \\
-0.063 & 0.027 & 0.078 & -0.005 & -0.056 & 0.019 \\
0.020 & -0.056 & -0.005 & 0.078 & 0.028 & -0.064 \\
0.019 & -0.005 & -0.056 & 0.028 & 0.078 & -0.064 \\
-0.039 & 0.020 & 0.019 & -0.064 & -0.064 & 0.127 \\
\end{psmallmatrix*},
\end{array}
$$
for a Gaussian random walk like in \eqref{Sigma_sRW_iid},
$$
\begin{array}{l@{\ }ll@{\ }l}
\text{\footnotesize true:} & \begin{psmallmatrix*}[r]
0.313 & -0.031 & -0.031 & -0.031 & -0.031 & -0.188 \\
-0.031 & 0.078 & 0.036 & -0.047 & -0.005 & -0.031 \\
-0.031 & 0.036 & 0.078 & -0.005 & -0.047 & -0.031 \\
-0.031 & -0.047 & -0.005 & 0.078 & 0.036 & -0.031 \\
-0.031 & -0.005 & -0.047 & 0.036 & 0.078 & -0.031 \\
-0.188 & -0.031 & -0.031 & -0.031 & -0.031 & 0.313 \\
\end{psmallmatrix*},
& \text{\footnotesize B348:} & 
\begin{psmallmatrix*}[r]
0.312 & -0.031 & -0.031 & -0.031 & -0.031 & -0.187 \\
-0.031 & 0.078 & 0.036 & -0.047 & -0.005 & -0.031 \\
-0.031 & 0.036 & 0.078 & -0.005 & -0.047 & -0.031 \\
-0.031 & -0.047 & -0.005 & 0.078 & 0.036 & -0.032 \\
-0.031 & -0.005 & -0.047 & 0.036 & 0.078 & -0.032 \\
-0.187 & -0.031 & -0.031 & -0.032 & -0.032 & 0.313 \\
\end{psmallmatrix*},
\\[-1ex]
\\
\text{\footnotesize U39:} & 
\begin{psmallmatrix*}[r]
0.313 & -0.031 & -0.031 & -0.031 & -0.031 & -0.188 \\
-0.031 & 0.078 & 0.036 & -0.047 & -0.005 & -0.031 \\
-0.031 & 0.036 & 0.078 & -0.005 & -0.047 & -0.031 \\
-0.031 & -0.047 & -0.005 & 0.078 & 0.036 & -0.031 \\
-0.031 & -0.005 & -0.047 & 0.036 & 0.078 & -0.031 \\
-0.188 & -0.031 & -0.031 & -0.031 & -0.031 & 0.313 \\
\end{psmallmatrix*},
& \text{\footnotesize U348:} & 
\begin{psmallmatrix*}[r]
0.313 & -0.032 & -0.031 & -0.031 & -0.031 & -0.188 \\
-0.032 & 0.078 & 0.036 & -0.047 & -0.005 & -0.031 \\
-0.031 & 0.036 & 0.078 & -0.005 & -0.047 & -0.031 \\
-0.031 & -0.047 & -0.005 & 0.078 & 0.036 & -0.032 \\
-0.031 & -0.005 & -0.047 & 0.036 & 0.078 & -0.032 \\
-0.188 & -0.031 & -0.031 & -0.032 & -0.032 & 0.314 \\
\end{psmallmatrix*},
\end{array}
$$
and for the coin-tossing DGP like in \eqref{Sigma_ct fair coin}:
$$
\begin{array}{l@{\ }ll@{\ }l}
\text{\footnotesize true:} & \begin{psmallmatrix*}[r]
0.313 & -0.031 & -0.031 & -0.031 & -0.031 & -0.188 \\
-0.031 & 0.078 & 0.016 & -0.047 & 0.016 & -0.031 \\
-0.031 & 0.016 & 0.078 & 0.016 & -0.047 & -0.031 \\
-0.031 & -0.047 & 0.016 & 0.078 & 0.016 & -0.031 \\
-0.031 & 0.016 & -0.047 & 0.016 & 0.078 & -0.031 \\
-0.188 & -0.031 & -0.031 & -0.031 & -0.031 & 0.313 \\
\end{psmallmatrix*},
& \text{\footnotesize B348:} & 
\begin{psmallmatrix*}[r]
0.314 & -0.031 & -0.031 & -0.032 & -0.031 & -0.188 \\
-0.031 & 0.078 & 0.016 & -0.047 & 0.015 & -0.031 \\
-0.031 & 0.016 & 0.078 & 0.015 & -0.047 & -0.031 \\
-0.032 & -0.047 & 0.015 & 0.078 & 0.016 & -0.031 \\
-0.031 & 0.015 & -0.047 & 0.016 & 0.078 & -0.031 \\
-0.188 & -0.031 & -0.031 & -0.031 & -0.031 & 0.313 \\
\end{psmallmatrix*},
\\[-1ex]
\\
\text{\footnotesize U39:} & 
\begin{psmallmatrix*}[r]
0.313 & -0.031 & -0.031 & -0.031 & -0.031 & -0.188 \\
-0.031 & 0.078 & 0.016 & -0.047 & 0.016 & -0.031 \\
-0.031 & 0.016 & 0.078 & 0.016 & -0.047 & -0.031 \\
-0.031 & -0.047 & 0.016 & 0.078 & 0.016 & -0.031 \\
-0.031 & 0.016 & -0.047 & 0.016 & 0.078 & -0.031 \\
-0.188 & -0.031 & -0.031 & -0.031 & -0.031 & 0.312 \\
\end{psmallmatrix*},
& \text{\footnotesize U348:} & 
\begin{psmallmatrix*}[r]
0.315 & -0.031 & -0.031 & -0.031 & -0.031 & -0.189 \\
-0.031 & 0.078 & 0.016 & -0.047 & 0.016 & -0.031 \\
-0.031 & 0.016 & 0.078 & 0.015 & -0.046 & -0.031 \\
-0.031 & -0.047 & 0.015 & 0.078 & 0.016 & -0.031 \\
-0.031 & 0.016 & -0.046 & 0.016 & 0.078 & -0.031 \\
-0.189 & -0.031 & -0.031 & -0.031 & -0.031 & 0.314 \\
\end{psmallmatrix*},
\end{array}
$$
We recognize a reasonable agreement between the true and the approximate covariance matrices in any case, where ``U39'' has the least computational effort among the three approximation methods.

%%%%%%%%%%%%%%%%

\section{Limit Theorems for Permutation Entropy and Statistical Complexity}
\label{section: statistics distribution}
Recall from Section~\ref{subsection: time series approach} that $\bfp = (p_1, \dots, p_{m!})^\top$ denotes the pmf-vector of the $m$th-order OPs, and that~$\hat{\bfp}$ are the corresponding relative frequencies as computed from the given OP-series $\Pi_1,\ldots,\Pi_n$. Furthermore, let~$\bfu = (1/m!, \ldots, 1/m!)^\top$ be the pmf vector of the discrete uniform distribution, and let~$\bfe_k$ refer to the $k$th one-point distribution, $k\in\{1,\ldots,m!\}$, i.e., with components $e_{k,l}=\delta_{kl}$. These vectors refer to the extreme scenarios of a pmf vector, namely perfect randomness versus perfect knowledge, and they also constitute the extremes of the entropy's range. Here, the (non-normalized) Shannon entropy of~$\bfp$, i.e., the permutation entropy according to \citet{bandt02}, is defined as
\begin{equation}
\label{eqEntropy}
H(\bfp)\ =\ 
-\sum_{k=1}^{m!} p_k\,\log{p_k}
\quad\text{with the convention}\quad
0\,\log{0} := 0,
\end{equation}
where the maximal value equals $H_0^{-1} := H(\bfu) = \log{m!}$, and where the minimal value~0 is attained for any $\bfp=\bfe_k$. With these notations, $H_0\cdot H(\bfp)$ equals the normalized permutation entropy.

\smallskip
In order to define the statistical complexity in the sense of \citet{rosso07}, let us next introduce the (non-normalized) Jensen--Shannon (JS) divergence (``disequilibrium'') as
\begin{eqnarray}
\label{eqJS}
D(\bfp) &=& 
H\big((\bfp+\bfu)/2\big) - \tfrac{1}{2}\,H(\bfp) - \tfrac{1}{2}\,H(\bfu)
\ =\ 
H\big((\bfp+\bfu)/2\big) - \tfrac{1}{2}\,H(\bfp) - \tfrac{1}{2}\,\log{m!}
\\
\nonumber
&&\text{with maximum}\quad
D_0^{-1} := D(\bfe_1) = -\frac{m!+1}{2 m!}\,\log(m!+1) + \log(2m!) - \tfrac{1}{2}\,\log{m!},
\end{eqnarray}
which expresses the distance of~$\bfp$ from the uniform distribution~$\bfu$. Finally, \citet{rosso07} define the (normalized) complexity as the product of the normalized entropy and JS-divergence:
\begin{equation}
\label{eqComplexity}
C(\bfp)\ =\ 
D_0\,D(\bfp)\cdot H_0\,H(\bfp).
\end{equation}
Then, we can summarize the dynamic structure (serial dependence) of the considered DGP $(X_t)_{t\in\bbz}$ by the point $\big(H_0\,H(\bfp),\ C(\bfp)\big)$ within the entropy--complexity plane. Recall from Section~\ref{intro} that the attainable region of $\big(H_0\,H(\bfp),\ C(\bfp)\big)$ is bounded by two curves (also see Figure~\ref{figure_HC_gauss_ct} below), which can be computed according to \citet{martin06} by using the R~package \href{https://cran.r-project.org/package=statcomp}{\nolinkurl{statcomp}}.

\medskip
In what follows, we are interested in the sample counterparts to the above quantities \eqref{eqEntropy}--\eqref{eqComplexity}, which are obtained by plugging-in the frequency vector~$\hat{\bfp}$ instead to the true OP-pmf~$\bfp$. From Section~\ref{section: OP frequencies distribution}, we know that for various types of DGP, the sample pmf $\hat{\bfp}$ is asymptotically normally distributed via $\sqrt{n} (\hat{\bfp}-\bfp) \xrightarrow{d} N(\mathbf{0}, \Sigma)$, where the covariance matrix~$\Sigma$ can be computed as shown in Section~\ref{section: autocovariance matrix}. Hence, our aim is to conclude from the asymptotics of~$\hat{\bfp}$ onto the (joint) asymptotics of $\big(H_0\,H(\hat{\bfp}),\ C(\hat{\bfp})\big)$. This, however, is more demanding as it might look like at first glance. As shown by \citet[Section~2]{wei_22}, it is crucial to not only focus on the first-order Taylor expansions of \eqref{eqEntropy}--\eqref{eqComplexity} as it is usually done when applying the ``Delta method'' \citep{serfling80}, but to consider the second-order Taylor expansions instead. More precisely, the second-order Taylor expansion of~$H\big(\hat{\bfp}\big)$ from \eqref{eqEntropy} around~$\bfp$ is given by
\begin{equation}
\label{PE_Taylor2}
H\big(\hat{\bfp}\big)\ \approx\ H(\bfp)
 - \sum_{k=1}^{m!} \big(1+\log{p_{k}}\big)\,\big(\hat{p}_k-p_{k}\big)
 - \sum_{k=1}^{m!} \frac{\big(\hat{p}_k-p_{k}\big)^2}{2\,p_{k}},
\end{equation}
where the linear term vanishes if $\bfp=\bfu$. Hence, the Delta method, which makes use of the first-order Taylor expansion
\begin{equation}
\label{PE_Taylor1}
H\big(\hat{\bfp}\big)\ \approx\ H(\bfp)
 - \sum_{k=1}^{m!} \big(1+\log{p_{k}}\big)\,\big(\hat{p}_k-p_{k}\big),
\end{equation}
is only applicable if $\bfp\not=\bfu$, which was also recognized by \citet[p.~3]{rey24}. This constitutes an important limitation as the uniform case $\bfp=\bfu$ happens for relevant DGP scenarios. First and foremost, if $(X_t)$ is i.i.d., then $\bfp=\bfu$ irrespective of the marginal distribution of~$X_t$, which allowed \citet{wei_22} to derive several non-parametric (distribution-free) tests for serial dependence. But there are also some dependence scenarios that lead to $\bfp=\bfu$, see \citet{bandt07} for a comprehensive discussion. This includes the MA$(q)$ process with equal weights and OP-length $m=2$, also recall our derivations in Section~\ref{subsection: MA(q)}, as well as certain stationary Gaussian processes if $m=3$. According to Propositions~4--5 in \citet{bandt07}, the OP-pmf~$\bfp$ is then of the form $(v, \frac{1-2v}{4}, \ldots, \frac{1-2v}{4}, v)$, where $v=\frac{1}{\pi}\,\arcsin\!\big(\frac{1}{2}\sqrt{\frac{1-\rho(2)}{1-\rho(1)}}\big)^\top$ with $\rho(h) = \Cor[X_t,X_{t-h}]$ denoting the autocorrelation function (acf). So $v=\frac{1}{6}$ and thus $\bfp=\bfu$ iff $\rho(2)=\rho(1)$. Hence, to sum up, when deriving the joint asymptotics of $\big(H_0\,H(\hat{\bfp}),\ C(\hat{\bfp})\big)$, we cannot exclude the case $\bfp=\bfu$ due to its practical relevance, but we have to consider both cases $\bfp\not=\bfu$ and $\bfp=\bfu$ separately.

\begin{figure}[t]
\centering\footnotesize
(a)\hspace{-3ex}\includegraphics[viewport=0 10 265 235, clip=, scale=0.75]{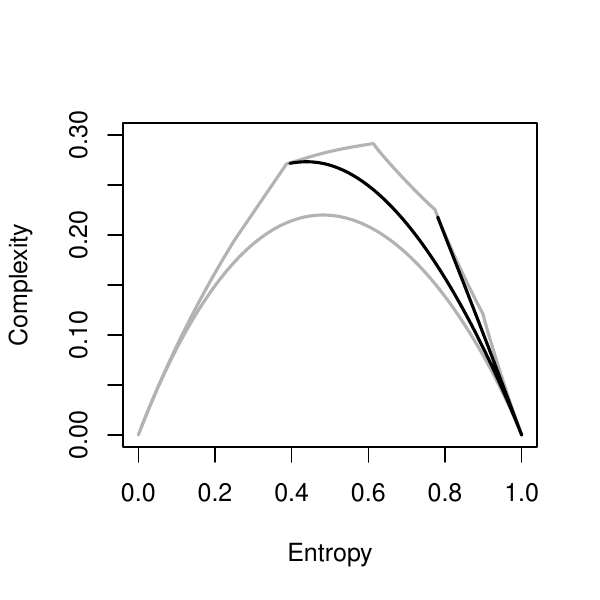}
\qquad
(b)\hspace{-3ex}\includegraphics[viewport=0 10 265 235, clip=, scale=0.75]{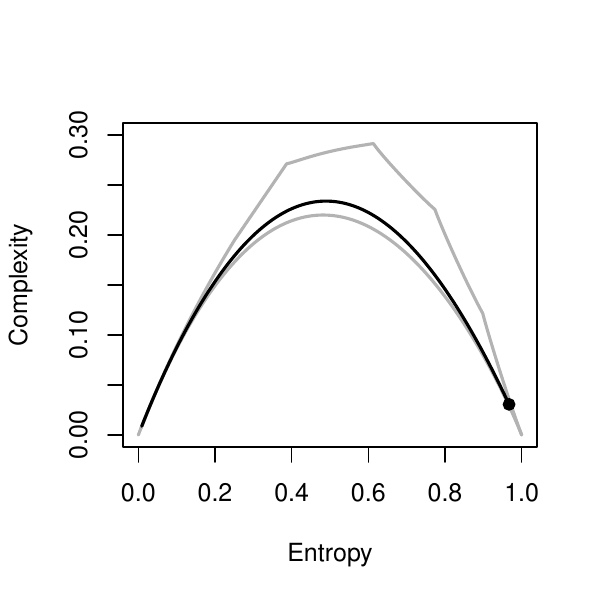}
\caption{Entropy--complexity plane for $m=3$, where bounds of attainable range shown in light grey. Trajectory of possible entropy--complexity pairs (black) for (a) stationary Gaussian processes and (b) GCT processes according to Proposition~\ref{prop: Sigma_ct}, where ``fair coin'' highlighted by point.}
\label{figure_HC_gauss_ct}
\end{figure}

\begin{remark}
\label{remark_HC_gauss_ct}
The possible entropy--complexity pairs for a stationary Gaussian process, i.e., with $\bfp = (v, \frac{1-2v}{4}, \ldots, \frac{1-2v}{4}, v)^\top$ for $v\in(0,0.5)$, are plotted in Figure~\ref{figure_HC_gauss_ct}\,(a). It is worth noting that this trajectory differs considerably from the one of a GCT process as plotted in part~(b), where $\bfp = (p^2, p^2 q, p^2 q, p q^2, p q^2, q^2)^\top$ with $p=1-q\in (0,1)$ according to Proposition~\ref{prop: Sigma_ct}. They only agree in one case, namely if $p=q=0.5$ (fair coin) and $v=0.25$.
\end{remark}

\subsection{Asymptotics for Non-Uniform OP~Distribution}
\label{Asymptotics for Non-Uniform OP-Distribution}
If $\bfp\not=\bfu$, then the Delta method \citep{serfling80} can be applied to the normal asymptotics
\begin{equation*}
\sqrt{n} (\hat{\bfp}-\bfp)\ \xrightarrow{d}\ N(\mathbf{0}, \Sigma)
\end{equation*}
according to Theorem~\ref{theorem: asymp distribution OP probs case 1} or~\ref{theorem: asymp distribution OP probs case 2}. While our ultimate aim are the joint asymptotics of $\big(H_0\,H(\hat{\bfp}),\ C(\hat{\bfp})\big)$, we start with the following auxiliary result.

\begin{lemma}
\label{lemmaH}
Let the assumptions of Theorem~\ref{theorem: asymp distribution OP probs case 1} or \ref{theorem: asymp distribution OP probs case 2} hold, and let $\bfp\not=\bfu$. Then, the entropy \eqref{eqEntropy} satisfies
\begin{equation}
\label{HH_joint_asymptotics}
\sqrt{n}\, \left(\begin{array}{c} H(\hat{\bfp}) - H(\bfp)\\ H\big((\hat{\bfp}+\bfu)/2\big) - H\big((\bfp+\bfu)/2\big) \end{array}\right)\ \xrightarrow{d}\ N(\mathbf{0}, \Sigma^{(1)}),
\end{equation}
where the entries of~$\Sigma^{(1)}$ are given by
\begin{eqnarray*}
\sigma^{(1)}_{11} &=& \textstyle
\sum_{i,j=1}^{m!} \log{p_i}\, \log{p_j}\, \sigma_{ij},
\\
\sigma^{(1)}_{12} &=& \textstyle
\frac{1}{2}\,\sum_{i,j=1}^{m!} \log{p_i}\, \log{\frac{p_j+1/m!}{2}}\, \sigma_{ij},
\\
\sigma^{(1)}_{22} &=& \textstyle
\frac{1}{4}\,\sum_{i,j=1}^{m!} \log{\frac{p_i+1/m!}{2}}\, \log{\frac{p_j+1/m!}{2}}\, \sigma_{ij}.
\end{eqnarray*}
Here, the $\sigma_{ij}$ are defined by \eqref{eq: long-run covariance matrix}.
\end{lemma}

The proof is given in Appendix~\ref{proof lemmaH}.
Note that the univariate normal asymptotics of solely $H(\hat{\bfp})$, as implied by Lemma~\ref{lemmaH}, have also been considered by \citet{rey24}. 
The covariance matrix $\Sigma^{(1)}$ has some interesting properties. First of all, it is singular for OPs of length $m=2$ regardless of the underlying distribution, as shown by the following lemma.

\begin{lemma}
\label{lemma: det Sigma1 m2}
    Let $m=2$. Then $\det \Sigma^{(1)} = 0$ holds true.
\end{lemma}

The proof is provided in Appendix~\ref{proof lemma: det Sigma1 m2}.
As a consequence of Lemma~\ref{lemma: det Sigma1 m2}, $H(\bfp)$ and $H\big((\bfp+\bfu)/2\big)$ are linearly dependent if $m=2$, which is reasonable as~$\bfp$ depends on only one parameter (recall that $p_2=1-p_1$). Consequently, $\bfp$ is already uniquely determined by $H(\bfp)$, up to a permutation of its components $p_1$ and $p_2$, which in turn gives $H(\bfp+\bfu)$.
Thus, for $m=2$, there is no information gain by considering the entropy--complexity pair instead of just the entropy, as we will see from the form of the long-run covariance in Theorem~\ref{theorem: asymp distr H C}. 

\smallskip
Now, the question arises as to whether $\det \Sigma^{(1)} = 0$ is also true for OPs of length $m \geq 3$. In general, the answer to this question is negative (although later in Section~\ref{section: applications}, we recognize only modest deviations from non-invertibility). A counterexample is the GCT process with $p=\tfrac{1}{4}$ and OPs of length $m=3$. From Proposition~\ref{prop: Sigma_ct}, we know that $\bfp=(\frac{1}{16},\frac{3}{64},\frac{3}{64},\frac{9}{64},\frac{9}{64},\frac{9}{16})^\top$, and that its long-run covariance matrix $\Sigma$ is given by
\begin{equation}
\label{ct025}
    \Sigma = \scalemath{0.9}{\begin{pmatrix}
        \tfrac{21}{4^4} &\tfrac{3}{4^5} &\tfrac{3}{4^5} &\tfrac{9}{4^5} &\tfrac{9}{4^5} &-\tfrac{27}{4^4} \\
        \tfrac{3}{4^5} &\tfrac{165}{4^6} &\tfrac{21}{4^6} &-\tfrac{81}{4^6} &\tfrac{63}{4^6} &-\tfrac{45}{4^5} \\
        \tfrac{3}{4^5} &\tfrac{21}{4^6} &\tfrac{165}{4^6} &\tfrac{63}{4^6} &-\tfrac{81}{4^6} &-\tfrac{45}{4^5} \\
        \tfrac{9}{4^5} &-\tfrac{81}{4^6} &\tfrac{63}{4^6} &\tfrac{333}{4^6} &\tfrac{189}{4^6} &-\tfrac{135}{4^5} \\
        \tfrac{9}{4^5} &\tfrac{63}{4^6} &-\tfrac{81}{4^6} &\tfrac{189}{4^6} &\tfrac{333}{4^6} &-\tfrac{135}{4^5} \\
        -\tfrac{27}{4^4} &-\tfrac{45}{4^5} &-\tfrac{45}{4^5} &-\tfrac{135}{4^5} &-\tfrac{135}{4^5} &\tfrac{117}{4^4}
    \end{pmatrix}}
    \quad\text{such that } \det{\Sigma^{(1)}}\approx 0.00184\ >0.
\end{equation}

\smallskip
Equation \eqref{HH_joint_asymptotics} in Lemma~\ref{lemmaH} allows to derive the joint asymptotic distribution of entropy~$H(\hat{\bfp})$ from \eqref{eqEntropy} and disequilibrium~$D(\hat{\bfp})$ from \eqref{eqJS}.

\begin{theorem}
\label{theorem: asymp distr H D}
Let the assumptions of Theorem~\ref{theorem: asymp distribution OP probs case 1} or \ref{theorem: asymp distribution OP probs case 2} hold, and let $\bfp\not=\bfu$. 
Then, the joint asymptotics of entropy \eqref{eqEntropy} and disequilibrium \eqref{eqJS} are
\begin{eqnarray*}
\sqrt{n}\, \left(\begin{array}{c} H(\hat{\bfp}) - H(\bfp)\\ D(\hat{\bfp}) - D(\bfp) \end{array}\right) &\xrightarrow{d}& N(\mathbf{0}, \Sigma^{(2)}),
\end{eqnarray*}
where the entries of~$\Sigma^{(2)}$ are given by
\begin{eqnarray*}
\sigma^{(2)}_{11} &=& \textstyle
\sum_{i,j=1}^{m!} \log{p_i}\, \log{p_j}\, \sigma_{ij},
\\
\sigma^{(2)}_{12} &=& \textstyle
%\frac{1}{2}\,\sum\limits_{i,j=1}^{m!} (1+ \log{p_i})(\log{\frac{p_j+1/m!}{2}}-\log{p_j})\, \sigma_{ij}
%\ =\ 
\frac{1}{2}\,\sum_{i,j=1}^{m!} \log{p_i}\,\log{\frac{p_j+1/m!}{2p_j}}\, \sigma_{ij},
\\
\sigma^{(2)}_{22} &=& \textstyle
%\frac{1}{4}\,\sum_{i,j=1}^{m!} \Big((1+ \log{\frac{p_i+1/m!}{2}})(1+\log{\frac{p_j+1/m!}{2}}) - (1+ \log{p_i})(1+\log{\frac{p_j+1/m!}{2}}) - (1+ \log{p_j})(1+\log{\frac{p_i+1/m!}{2}}) + (1+ \log{p_i})(1+\log{p_j})\Big)\, \sigma_{ij}.
%\frac{1}{4}\,\sum\limits_{i,j=1}^{m!} (\log{\frac{p_i+1/m!}{2}}-\log{p_i})(\log{\frac{p_j+1/m!}{2}}-\log{p_j})\, \sigma_{ij}.
\frac{1}{4}\,\sum_{i,j=1}^{m!} \log{\frac{p_i+1/m!}{2 p_i}}\, \log{\frac{p_j+1/m!}{2 p_j}}\, \sigma_{ij}.
\end{eqnarray*}
Here, the $\sigma_{ij}$ are defined by \eqref{eq: long-run covariance matrix}.
\end{theorem}

The proof is provided in Appendix~\ref{proof theorem: asymp distr H D}. 
Finally, we compute the joint asymptotic distribution of normalized entropy \eqref{eqEntropy} and complexity \eqref{eqComplexity}, as these are plotted against each other in the entropy--complexity plane. 

\begin{theorem}
\label{theorem: asymp distr H C}
Let the assumptions of Theorem~\ref{theorem: asymp distribution OP probs case 1} or \ref{theorem: asymp distribution OP probs case 2} hold, and let $\bfp\not=\bfu$. 
Then, the joint asymptotics of normalized entropy \eqref{eqEntropy} and complexity \eqref{eqComplexity} are
\begin{eqnarray*}
\sqrt{n}\, \left(\begin{array}{c} H_0\,\big(H(\hat{\bfp}) - H(\bfp)\big)\\ C(\hat{\bfp}) - C(\bfp) \end{array}\right)\ \xrightarrow{d}\ N(\mathbf{0}, \Sigma^{(3)}),
\end{eqnarray*}
where the entries of~$\Sigma^{(3)}$ are given by
\begin{eqnarray*}
\sigma^{(3)}_{11} &=& 
%H_0^2\,\sigma^{(1)}_{11}
%\ =\ 
H_0^2\,\sigma^{(2)}_{11},
\\
\sigma^{(3)}_{12} &=& 
%H_0^2D_0\,\big((D(\bfp) - \tfrac{1}{2}\,H(\bfp))\,\sigma^{(1)}_{11} + H(\bfp)\,\sigma^{(1)}_{12}\big)
%\ =\ H_0^2D_0\,\big(D(\bfp)\,\sigma^{(1)}_{11} + H(\bfp)\,(\sigma^{(1)}_{12} - \tfrac{1}{2}\,\sigma^{(1)}_{11})\big),
%\ =\ 
H_0^2D_0\,\big(D(\bfp)\,\sigma^{(2)}_{11} + H(\bfp)\,\sigma^{(2)}_{12}\big),
\\
\sigma^{(3)}_{22} &=& 
%H_0^2D_0^2\,\big((D(\bfp) - \tfrac{1}{2}\,H(\bfp))^2\,\sigma^{(1)}_{11} + 2\,H(\bfp)\,(D(\bfp) - \tfrac{1}{2}\,H(\bfp))\,\sigma^{(1)}_{12} + H(\bfp)^2\,\sigma^{(1)}_{22}\big)
%\\
 %&=& 
H_0^2D_0^2\,\big(D(\bfp)^2\,\sigma^{(2)}_{11} + 2\,H(\bfp)\,D(\bfp)\,\sigma^{(2)}_{12} + H(\bfp)^2\,\sigma^{(2)}_{22}\big).
\end{eqnarray*}
Here, the $\sigma^{(2)}_{ij}$ are provided by Theorem~\ref{theorem: asymp distr H D}.
\end{theorem}

The proof is given in Appendix~\ref{proof theorem: asymp distr H C}. Recall that our asymptotics in Theorem~\ref{theorem: asymp distr H C} differ from those derived by \citet{rey25}, because the latter authors did not consider the true OP-process but a ``multinomial model'' instead.

\subsection{Asymptotics for Uniform OP~Distribution}
\label{Asymptotics for Uniform OP-Distribution}
Next let us turn to the case $\bfp=\bfu$, where the first-order Taylor expansion \eqref{PE_Taylor1} degenerates to the constant expression $H(\bfu)=\log{m!}$. Therefore, we make use of the second-order Taylor expansion \eqref{PE_Taylor2} for our asymptotic derivations, which simplifies for $\bfp=\bfu$ to
\begin{equation}
\label{PE_Taylor2_uniform}
\textstyle
H\big(\hat{\bfp}\big)\ \approx\ \log{m!}
 - \frac{m!}{2}\,\sum_{k=1}^{m!} \big(\hat{p}_k-p_{k}\big)^2,
\end{equation}
also see \citet{wei_22}. Together with the normal asymptotics from Theorem~\ref{theorem: asymp distribution OP probs case 1} or \ref{theorem: asymp distribution OP probs case 2}, it follows that
\begin{equation}
\label{H_QFasymp}
\textstyle
n\, (-\tfrac{2}{m!})\,\big(H\big(\hat{\bfp}\big) - \log{m!}\big)\ \xrightarrow{d}\ \sum_{i=1}^l \lambda_i\,\chi_{r_i}^2\ =:Q_m,
\end{equation}
where $\lambda_1,\ldots,\lambda_l$ are the distinct non-zero eigenvalues of~$\Sigma$ with algebraic multiplicities $r_1,\ldots,r_l$, where~$\chi_r^2$ refers to a $\chi^2$-distributed random variable with $r$ degrees of freedom, and where $\chi_{r_1}^2,\ldots,\chi_{r_l}^2$ are independent of each other. This argumentation from \citet{wei_22} is based on Theorem~3.1 in \citet{tan_77}. The distribution of~$Q_m$, which is referred to as a quadratic-form (QF) distribution, can be evaluated numerically by using the R~package \href{https://cran.r-project.org/package=CompQuadForm}{\nolinkurl{CompQuadForm}} of \citet{duchesne10}.

\smallskip
Next, we adapt the above argumentation to the statistics~$D(\hat{\bfp})$ from \eqref{eqJS} as well as~$C(\hat{\bfp})$ from \eqref{eqComplexity}. 

\begin{theorem}
\label{theorem: asymp distr H D uniform}
Let the assumptions of Theorem~\ref{theorem: asymp distribution OP probs case 1} or \ref{theorem: asymp distribution OP probs case 2} hold, and let $\bfp=\bfu$. 
Then, the joint asymptotics of entropy \eqref{eqEntropy} and disequilibrium \eqref{eqJS} are
\begin{equation*}
%\label{HD_joint_asymptotics_uniform}
n\, \left(\begin{array}{c} H(\hat{\bfp}) - \log{m!}\\ D(\hat{\bfp}) \end{array}\right)\ \xrightarrow{d}\ \left(\begin{array}{c} -\frac{m!}{2}\\ \frac{m!}{8} \end{array}\right)\, Q_m,
\end{equation*}
where~$Q_m$ is the QF-distributed random variable from \eqref{H_QFasymp}.
\end{theorem}

The proof is provided in Appendix~\ref{proof theorem: asymp distr H D uniform}.
Note that Theorem~\ref{theorem: asymp distr H D uniform} implies an asymptotic linear relationship between~$H(\hat{\bfp})$ and~$D(\hat{\bfp})$, where the absolute deviation from the point $(\log{m!},0)$ follows a rescaled QF-distribution. 

\begin{theorem}
\label{theorem: asymp distr H C uniform}
Let the assumptions of Theorem~\ref{theorem: asymp distribution OP probs case 1} or \ref{theorem: asymp distribution OP probs case 2} hold, and let $\bfp=\bfu$. 
Then, the joint asymptotics of normalized entropy \eqref{eqEntropy} and complexity \eqref{eqComplexity} are
\begin{eqnarray*}
%\label{HC_joint_asymptotics_uniform}
n\, \left(\begin{array}{c} H(\hat{\bfp})/\log{m!} - 1\\ C(\hat{\bfp}) \end{array}\right)\ \xrightarrow{d}\ \left(\begin{array}{c} -\frac{m!}{2\,\log{m!}}\\ \frac{m!}{8}\, D_0 \end{array}\right)\, Q_m,
\end{eqnarray*}
where~$Q_m$ is the QF-distributed random variable from \eqref{H_QFasymp}.
\end{theorem}

The proof is provided in Appendix~\ref{proof theorem: asymp distr H C uniform}.
Theorem~\ref{theorem: asymp distr H C uniform} implies that asymptotically, the estimated normalized entropy and complexity are on the straight line $f(H)=\frac{\log{m!}}{4}\, D_0\, (1-H)$ within the entropy--complexity plane, where the absolute deviations from the point $(1,0)$ follows a rescaled QF-distribution.

%%%%%%%%%%%%%%%%%%%%%%%%%%%%%%%%%%%%%%%%%%%%%%%%%%%%

\section{Applications}
\label{section: applications}

\begin{figure}[t]
\centering\footnotesize
(a)\hspace{-3ex}\includegraphics[viewport=0 10 265 230, clip=, scale=0.75]{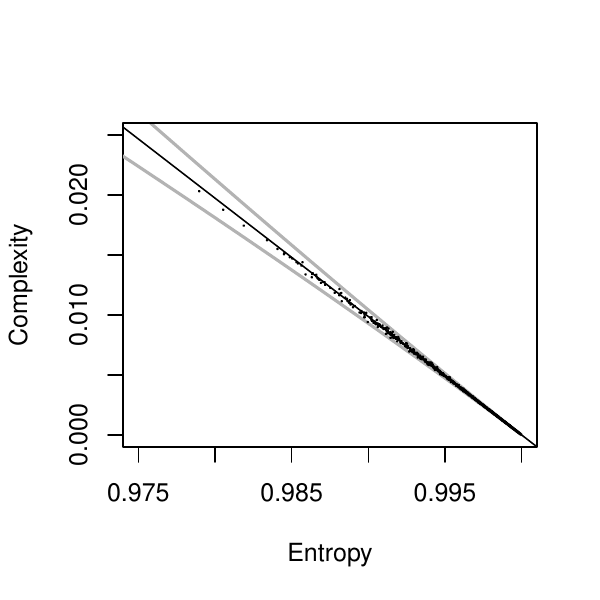}
\qquad
(b)\hspace{-3ex}\includegraphics[viewport=0 10 265 230, clip=, scale=0.75]{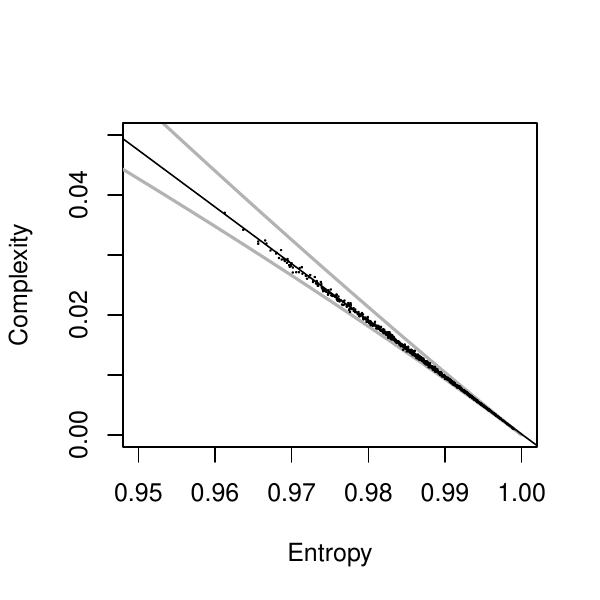}
\\[2ex]
(c)\hspace{-3ex}\includegraphics[viewport=0 10 265 230, clip=, scale=0.75]{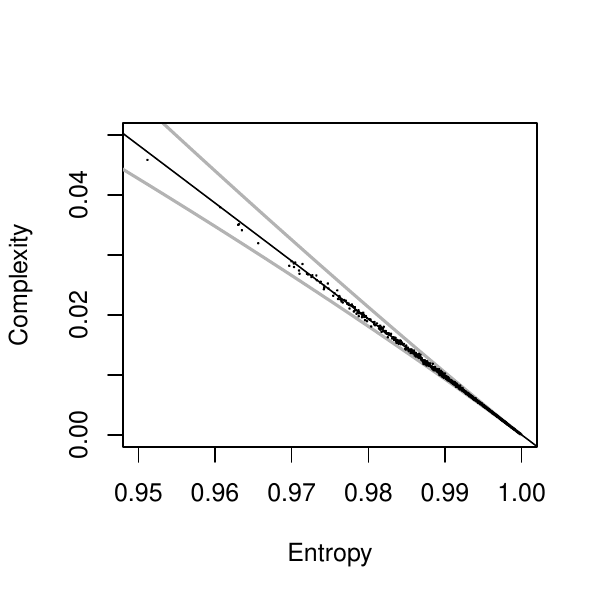}
\qquad
(d)\hspace{-3ex}\includegraphics[viewport=0 10 265 230, clip=, scale=0.75]{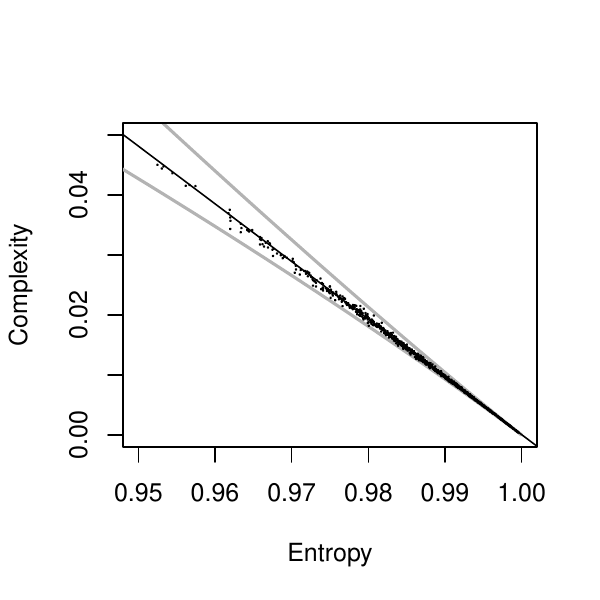}
\caption{Detail of entropy--complexity plane ($m=3$), where bounds of attainable range shown in light grey. Plot of 1,000 simulated statistics (sample size $T=250$) of normalized entropy \eqref{eqEntropy} and complexity \eqref{eqComplexity}. (a) i.i.d.\ DGP compared to line $\frac{\log{m!}}{4}\, D_0\, \big(1-H_0\, H(\hat\bfp)\big)$. (b) AR$(1)$ DGP, (c) QMA$(1)$ DGP, and (d) TEAR$(1)$ DGP, compared to respective line $C(\bfp)+\sigma_{12}^{(3)}/\sigma_{11}^{(3)}\, H_0\big(H(\hat\bfp)-H(\bfp)\big)$.}
\label{figure_plot_HCplane_tests}
\end{figure}

\subsection{Testing for Serial Dependence}
\label{subsection: testing dependence}
An obvious application of our newly derived asymptotic distributions in Section~\ref{section: statistics distribution} is their use for testing for serial dependence in the given time-series data. More precisely, the aim is to test the null hypothesis of serial independence against the alternative of a serially dependent DGP. Recall that under the i.i.d.-null, the OP-series always has a marginal uniform distribution, so we are in the case $\bfp=\bfu$ discussed in Section~\ref{Asymptotics for Uniform OP-Distribution}. In Theorems~\ref{theorem: asymp distr H D uniform} and~\ref{theorem: asymp distr H C uniform}, we obtained the notable result that both pairs, $(H,D)$ and $(H,C)$, are asymptotically distributed along a straight line, where the deviation from the null value $(\log{m!},0)$ or $(1,0)$, respectively, follows a rescaled QF-distribution. This strong linear relation between entropy and complexity is also illustrated by Figure~\ref{figure_plot_HCplane_tests}\,(a), where 1,000 simulated statistics of normalized entropy \eqref{eqEntropy} and complexity \eqref{eqComplexity} with $m=3$ are compared to the straight line $\frac{\log{m!}}{4}\, D_0\, \big(1-H_0\, H(\hat\bfp)\big)$ derived from Theorem~\ref{theorem: asymp distr H C uniform}. We recognize only modest deviations from this straight line, so the true sample distribution is well described by the linear asymptotic relationship according to Theorems~\ref{theorem: asymp distr H C uniform}. 
Therefore, it is reasonable to transform the vectors $(H,D)$ and $(H,C)$ into a scalar statistic in an appropriate way. The probably most obvious approach is a linear transformation like
\begin{equation}
\label{HDtest}
\widehat{HD}\ :=\ \frac{n}{m!}\,\Big(\log{m!}-H(\hat{\bfp})+4\, D(\hat{\bfp})\Big)
\end{equation}
concerning Theorem~\ref{theorem: asymp distr H D uniform}, and
\begin{equation}
\label{HCtest}
\widehat{HC}\ :=\ \frac{n}{m!}\,\Big(\log{m!}-H(\hat{\bfp})+\frac{4}{D_0}\, C(\hat{\bfp})\Big)
\end{equation}
concerning Theorem~\ref{theorem: asymp distr H C uniform}, which both follow the QF-distribution from \eqref{H_QFasymp} under the i.i.d.-null. We also constructed test statistics based on the taxi norm and Euclidean norm of $(H,D)$ and $(H,C)$, respectively, but the obtained size and power values are virtually identical to those of statistics \eqref{HDtest} and \eqref{HCtest}. Therefore, we omit these statistics for the subsequent analyses.

\begin{table}[t]
\centering
\caption{Rejection rates of entropy-disequilibrium test (``$HD$'') and entropy-complexity test (``$HC$'') for $m=3$, with 10,000 replications per DGP. Rejection rates of entropy test (``$H$'') taken from \citet{wei_22} and written in italic font.}
\label{tabTests}

\smallskip
\resizebox{\linewidth}{!}{
\begin{tabular}{l|ccc|ccc|ccc|ccc}
\toprule
DGP & \multicolumn{3}{c|}{i.i.d.} & \multicolumn{3}{c|}{AR$(1)$} & \multicolumn{3}{c|}{QMA$(1)$} & \multicolumn{3}{c}{TEAR$(1)$} \\
$T$ & $HD$ & $HC$ & $H$ & $HD$ & $HC$ & $H$ & $HD$ & $HC$ & $H$ & $HD$ & $HC$ & $H$ \\
\midrule
100 & 0.055 & 0.053 & \it 0.051 & 0.256 & 0.252 & \it 0.250 & 0.160 & 0.156 & \it 0.157 & 0.248 & 0.244 & \it 0.246 \\
250 & 0.049 & 0.048 & \it 0.054 & 0.615 & 0.611 & \it 0.614 & 0.343 & 0.341 & \it 0.330 & 0.527 & 0.524 & \it 0.534 \\
500 & 0.051 & 0.051 & \it 0.051 & 0.933 & 0.933 & \it 0.929 & 0.604 & 0.601 & \it 0.587 & 0.835 & 0.833 & \it 0.824 \\
1000 & 0.050 & 0.050 & \it 0.049 & 1.000 & 1.000 & \it 0.999 & 0.886 & 0.886 & \it 0.888 & 0.988 & 0.988 & \it 0.984 \\
\bottomrule
\end{tabular}}
\end{table}

\smallskip
In order to investigate the finite-sample performance of the tests \eqref{HDtest} and \eqref{HCtest}, we did a simulation study with 10,000 replications per scenario. Based on OPs of length $m=3$, we always determined the empirical rejection rate, which expresses the size for a null DGP, and the power for an alternative DGP. Here, we used DGPs from the study of \citet{wei_22} such that our newly obtained size and power values are directly comparable to those of \citet{wei_22}, see Table~\ref{tabTests} for a summary. More precisely, the null DGP is i.i.d.\ $N(0,1)$, whereas our three alternative DGPs are
\begin{itemize}
    \item the first-order autoregressive (AR$(1)$) process $X_t=0.5\,X_{t-1}+\varepsilon_t$ with i.i.d.\ $\varepsilon_t\sim N(0,1)$;
    \item the first-order quadratic moving-average (QMA$(1)$) process $X_t=\varepsilon_t + 0.8\,\varepsilon_{t-1}^2$ with i.i.d.\ $\varepsilon_t\sim N(0,1)$;
    \item the TEAR$(1)$ process $X_t=B_t\,X_{t-1}+0.85\,\varepsilon_t$ with i.i.d.\ $\varepsilon_t\sim \text{Exp}(1)$ and Bernoulli $B_t$ with $P(B_t=1)=0.15$.
\end{itemize}
Further details and characteristics about these DGPs can be found in \citet{wei_22}. But these are not required for the interpretation of Table~\ref{tabTests}, because its main message is as follows: There is hardly any difference in size and power if using~$\widehat{HD}$ or~$\widehat{HC}$ instead of~$H(\hat{\bfp})$ (beyond simulation uncertainty). The additional information provided by~$D(\hat{\bfp})$ or~$C(\hat{\bfp})$ compared to solely~$H(\hat{\bfp})$ does not help in uncovering the alternative scenario, the relevant information is already provided by the entropy itself. This is demonstrated by parts~(b)--(d) of Figure~\ref{figure_plot_HCplane_tests}, where 1,000 simulated statistics for each alternative DGP are plotted in the entropy--complexity plane. Parts~(b)--(d) compare the pairs to the straight line being implied by the normal asymptotics of Theorem~\ref{theorem: asymp distr H C}, namely $C(\bfp)+\sigma_{12}^{(3)}/\sigma_{11}^{(3)}\, H_0\big(H(\hat\bfp)-H(\bfp)\big)$. Here, the actual expressions for~$\bfp$, $H_0\, H(\bfp)$, $C(\bfp)$, and~$\Sigma^{(3)}$ are
\begin{itemize}
    \item AR$(1)$:\quad $\bfp\approx (0.210, 0.145, 0.145, 0.145, 0.145, 0.210)^\top$, $H_0\, H(\bfp)\approx 0.9910$, $C(\bfp)\approx 0.0087$, and $\Sigma^{(3)}\approx \begin{psmallmatrix} 0.0087 & -0.0083 \\ -0.0083 & 0.0078 \end{psmallmatrix}$;
    \item QMA$(1)$:\quad $\bfp\approx (0.150, 0.170, 0.180, 0.151, 0.142, 0.207)^\top$, $H_0\, H(\bfp)\approx 0.9951$, $C(\bfp)\approx 0.0048$, and $\Sigma^{(3)}\approx \begin{psmallmatrix} 0.0084 & -0.0081 \\ -0.0081 & 0.0078\end{psmallmatrix}$;
    \item TEAR$(1)$:\quad $\bfp\approx (0.218, 0.137, 0.148, 0.185, 0.174, 0.137)^\top$, $H_0\, H(\bfp)\approx 0.9916$, $C(\bfp)\approx 0.0082$, and $\Sigma^{(3)}\approx \begin{psmallmatrix} 0.0152 & -0.0146 \\ -0.0146 & 0.0140 \end{psmallmatrix}$;
\end{itemize}
which were computed by using the simulation-based covariance approximation of Section~\ref{Approximate Computation of Covariance Matrix}. It can be seen that also under the alternative, the statistics exhibit a nearly perfect linear relationship such that the value of~$C(\hat\bfp)$ is almost perfectly implied by the value of~$H(\hat\bfp)$.

\subsection{Estimation Uncertainty in Entropy--Complexity Plane}
\label{Estimation Uncertainty in Entropy--Complexity Plane}
\begin{figure}[t]
\centering\footnotesize
(a)\hspace{-3ex}\includegraphics[viewport=0 10 265 230, clip=, scale=0.75]{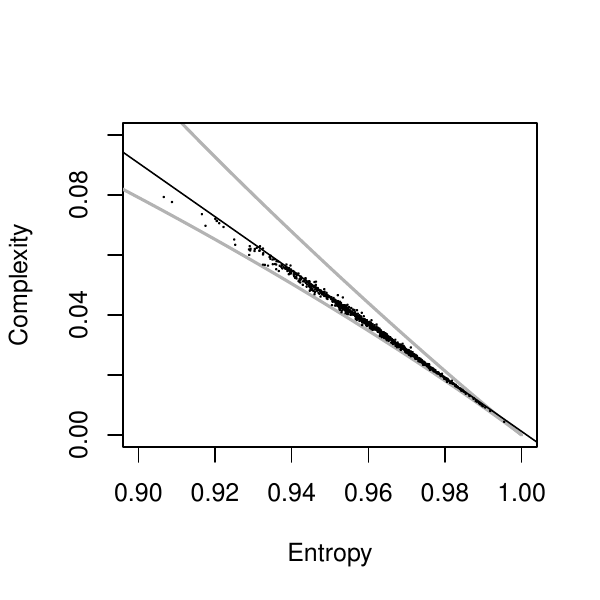}
\qquad
(b)\hspace{-3ex}\includegraphics[viewport=0 10 265 230, clip=, scale=0.75]{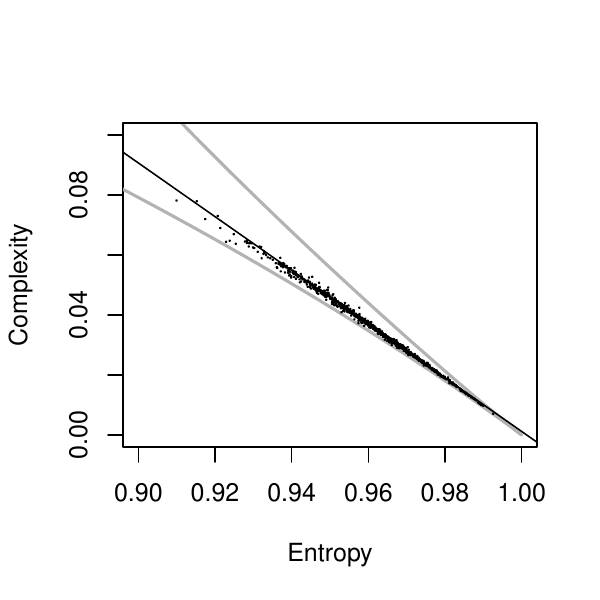}
\caption{Detail of entropy--complexity plane ($m=3$), where bounds of attainable range shown in light grey. Plot of 1,000 simulated statistics (sample size $T=250$) of normalized entropy \eqref{eqEntropy} and complexity \eqref{eqComplexity}. (a) Gaussian random walk and (b) coin-tossing process, compared to common line $C(\bfp)+\sigma_{12}^{(3)}/\sigma_{11}^{(3)}\, H_0 \big(H(\hat\bfp)-H(\bfp)\big) \approx 0.896 - 0.895\,H_0\, H(\hat\bfp)$.}
\label{figure_plot_HCplane_uncertainty}
\end{figure}

The normal asymptotics of Section~\ref{Asymptotics for Non-Uniform OP-Distribution}, especially Theorem~\ref{theorem: asymp distr H C}, can be used to approximate the estimation uncertainty of the pairs $\big(H(\hat\bfp), C(\hat\bfp)\big)$. As our first illustrative example, we consider a Gaussian random walk on the one hand, and the ``fair'' coin-tossing process on the other hand. While both lead to the same marginal distribution of the OP-series with $m=3$, namely $\bfp=(\frac{1}{4}, \frac{1}{8}, \ldots, \frac{1}{8}, \frac{1}{4})^\top$ with $H(\bfp)\approx 0.9671$ and $C(\bfp)\approx 0.0306$, they have slightly different covariance matrices~$\Sigma$, recall \eqref{Sigma_sRW_iid} versus \eqref{Sigma_ct fair coin}. Interestingly, however, if we compute the covariance matrices $\Sigma^{(1)}$,  $\Sigma^{(2)}$, and $\Sigma^{(3)}$ from Section~\ref{Asymptotics for Non-Uniform OP-Distribution}, they are identical and given by
\begin{equation}
\label{ct_grw_covariances}
\Sigma^{(1)} \approx \begin{psmallmatrix} 0.1201 & 0.0309 \\ 0.0309 & 0.0080 \end{psmallmatrix},
\qquad
\Sigma^{(2)} \approx \begin{psmallmatrix} 0.1201 & -0.0292 \\ -0.0292 & 0.0071 \end{psmallmatrix},
\qquad
\Sigma^{(3)} \approx \begin{psmallmatrix} 0.0374 & -0.0335 \\ -0.0335 & 0.0300 \\ \end{psmallmatrix}.
\end{equation}
$\Sigma^{(1)}$, $\Sigma^{(2)}$, and $\Sigma^{(3)}$ are not invertible, with the truly positive eigenvalue being given by about 0.1281, 0.1272, and 0.0674, respectively. The zero eigenvalue is caused by the asymptotically linear relationship between the pairs of statistics, also recall our previous discussion. In analogy to Section~\ref{subsection: testing dependence}, this linear relationship is also confirmed empirically, see Figure~\ref{figure_plot_HCplane_uncertainty}, where we recognize only slight deviations from the common line $C(\bfp)+\sigma_{12}^{(3)}/\sigma_{11}^{(3)}\, H_0\big(H(\hat\bfp)-H(\bfp)\big) \approx 0.896 - 0.895\,H_0\,H(\hat\bfp)$.

\begin{figure}[t]
\centering\footnotesize
(a)\hspace{-3ex}\includegraphics[viewport=0 10 265 230, clip=, scale=0.75]{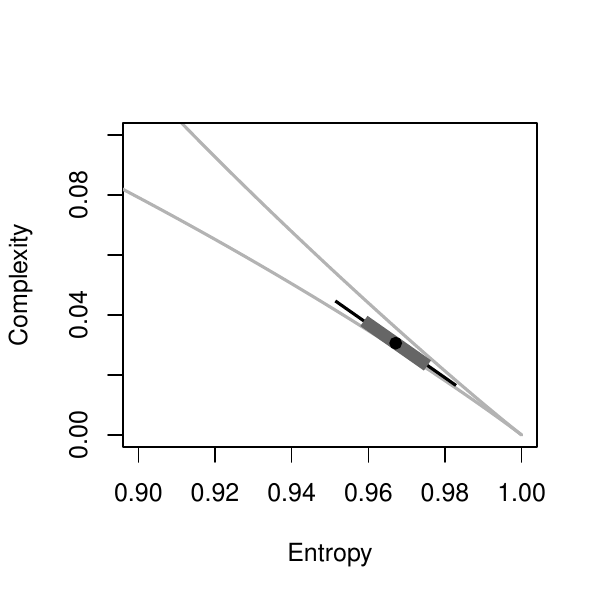}
\qquad
(b)\hspace{-3ex}\includegraphics[viewport=0 10 265 230, clip=, scale=0.75]{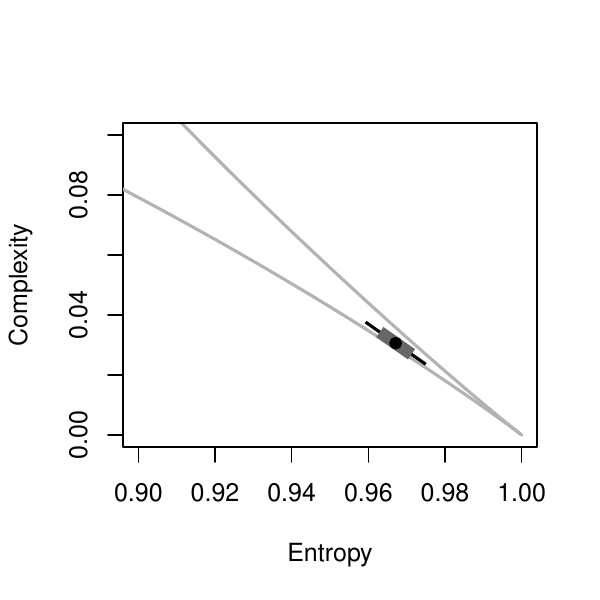}
\caption{Detail of entropy--complexity plane ($m=3$), where bounds of attainable range shown in light grey. Kind of boxplot for expressing estimation uncertainty, built up of asymptotic deciles, quartiles, and median for Gaussian random walk or coin-tossing process. Sample size (a) $T=250$ and (b) $T=1000$.}
\label{figure_plot_quantiles}
\end{figure}

\smallskip
The normal asymptotics \eqref{ct_grw_covariances} can be utilized in various ways to quantify the estimation uncertainty of entropy and complexity. Approximate standard errors immediately follow as $\sqrt{0.0374/n}$ for $H_0\,H(\hat\bfp)$ and $\sqrt{0.0300/n}$ for $C(\hat\bfp)$. Furthermore, as a visual solution for expressing the estimation uncertainty (and considering that~$\Sigma^{(3)}$ is not invertible), we suggest to construct a kind of boxplot, see Figure~\ref{figure_plot_quantiles}. The asymptotic median is represented by a black dot, the thick grey line is plotted between the quartiles, and the thin black line between the deciles. From this graph, the orientation of the estimated entropy--complexity pairs immediately gets clear, and also the effect of the sample size~$T$ can be recognized by comparing parts~(a) and~(b). 

\begin{figure}[t]
\centering\footnotesize
(a)\hspace{-3ex}\includegraphics[viewport=0 10 265 235, clip=, scale=0.75]{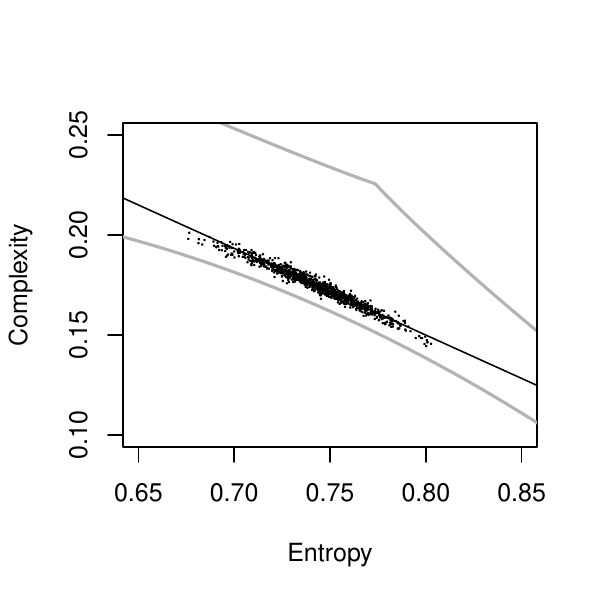}
\qquad
(b)\hspace{-3ex}\includegraphics[viewport=0 10 265 235, clip=, scale=0.75]{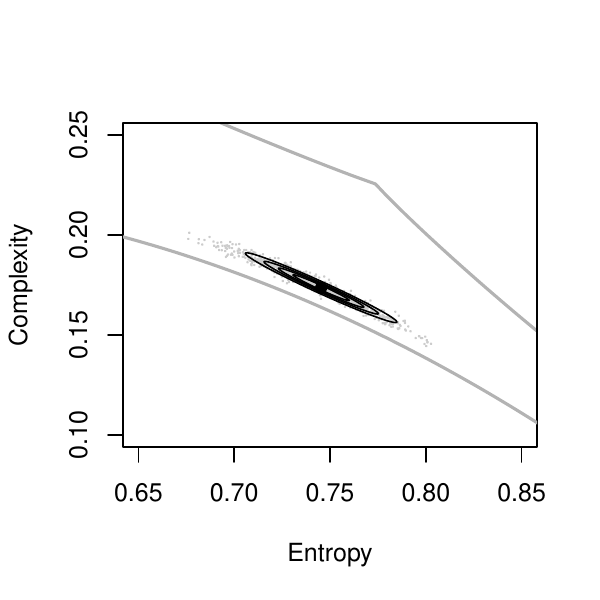}
\caption{Detail of entropy--complexity plane ($m=3$), where bounds of attainable range shown in light grey. For coin-tossing process with $p=0.25$, (a) shows plot of 1,000 simulated statistics ($n=1000$) of normalized entropy \eqref{eqEntropy} and complexity \eqref{eqComplexity} compared to line $C(\bfp)+\sigma_{12}^{(3)}/\sigma_{11}^{(3)}\, H_0 \big(H(\hat\bfp)-H(\bfp)\big) \approx 0.497 - 0.434\,H_0\, H(\hat\bfp)$, and (b) corresponding ellipsoids derived from bivariate normal distribution according to Theorem~\ref{theorem: asymp distr H C} and Equation \eqref{ct025}.}
\label{figure_ct025}
\end{figure}

\smallskip
The above illustrative examples led to a non-invertible~$\Sigma^{(3)}$, for which reason we used a kind of boxplot for visualizing the estimation uncertainty. But as shown in Section~\ref{Asymptotics for Non-Uniform OP-Distribution}, if the OPs' length is $m=3$, then it is also possible that~$\Sigma^{(3)}$ is truly positive definite, which was shown by the GCT process with $p=0.25$, recall \eqref{ct025}. Therefore, in such invertible cases, we can visualize the estimation uncertainty of the pairs $\big(H(\hat\bfp), C(\hat\bfp)\big)$ in a more traditional way, namely by computing ellipsoids from the truly bivariate density function of their asymptotic normal distribution. This is illustrated by Figure~\ref{figure_ct025}, where the ellipsoids fit quite well to the 1,000 simulated statistics. It should be noted, however, that although we are concerned with a truly bivariate normal distribution this time, the deviations from the straight line $C(\bfp)+\sigma_{12}^{(3)}/\sigma_{11}^{(3)}\, H_0 \big(H(\hat\bfp)-H(\bfp)\big) \approx 0.497 - 0.434\,H_0\, H(\hat\bfp)$ are small. In fact, the asymptotic cross-correlation between normalized entropy \eqref{eqEntropy} and complexity \eqref{eqComplexity} is given by $\approx -0.983$ and thus still very close to perfect negative correlation. Hence, the additional information provided by~$C(\hat{\bfp})$ compared to solely~$H(\hat{\bfp})$ is modest also in this invertible case.

\section{Conclusions}
\label{Conclusions}
In this research, we accepted the challenge of \citet[p.~394]{chagas22}, ``Making statistical inference in the entropy–complexity plane is challenging, but its potential impact is huge.'', and derived closed-form joint asymptotics for (among others) normalized entropy and complexity. To make this possible, we first derived some limit theorems for OP~frequencies under weak dependence conditions on the underlying process, where we distinguished between two perspectives for our limit theorems, namely the time-series and the OP-series perspective. For both cases, we derived novel closed-form expressions for the OP~frequencies' long-run covariance, and we also provided a simulation-based approach for approximate covariance computation. In context of the OP-series perspective and extending an earlier contribution by \citet{bandt25}, we proposed the GCT process and derived general closed-form asymptotics for the resulting OP~frequencies. Regarding the joint asymptotics for normalized entropy and complexity, we distinguished the case of a non-uniform OP distribution where we obtained normal asymptotics, and that of a uniform OP distribution where we obtained QF-asymptotics. It became clear, also when considering some applications, that normalized entropy and complexity asymptotically exhibit very strong up to perfect negative correlation. For this reason, dependence tests based on both normalized entropy and complexity turned out to not be more powerful than tests based on solely the entropy. Finally, we applied the novel asymptotics for (visually) expressing the estimation uncertainty in the entropy--complexity plane.

\smallskip
There are various directions for future research. First, it would be relevant to extend our asymptotic derivations to random fields and the corresponding spatial OPs, where only the i.i.d.\ case has been considered so far \citep[see][]{wei_kim_24}. Possible applications of spatial OPs together with the entropy--complexity plane are described by \citet{ribeiro12,zunino16,sigaki18}, for example. Second, the novel GCT process deserves further attention. It constitutes one of the few model families, where closed-form asymptotics are fully available, and which allows to reach many regions within the entropy--complexity plane. The GCT process could serve as the base for deriving the asymptotics of further OP-statistics, and it could also be investigated for OPs of length $m\geq 4$.

\subsubsection*{Acknowledgments}
%The authors thank the two referees for their useful comments on an earlier draft of this article.
Both authors are grateful to Professor Annika Betken (University of Twente) and Professor Alexander Schnurr (University of Siegen) for inviting them to the COPAB workshop, in the course of which the presented research was inspired.

%\newpage

%\nocite{*}
\bibliography{references}
\bibliographystyle{apalike}

\appendix
%%%%%%%%%%%%%%%%%%%%%%%
\section{Approximating Functionals and Mixing Conditions}
\label{section: approximating functionals and mixing conditions}

In this section, we give a short overview on concepts of short-range dependence employed in this paper. We begin with two mixing conditions.

\begin{definition}
    Let $(Z_t)_{t\in\bbz}$ be a process defined on a probability space $(\Omega, \mathcal{F}, \bbp)$. For $k \leq l$, define the $\sigma$-fields $\ca_{k}^l:=\sigma(Z_k,\ldots,Z_l)$.
    \begin{enumerate}[label=(\alph*)]
        \item $(Z_t)_{t\in\bbz}$ is called \emph{absolutely regular} if $\beta_k\to 0$ where
        \begin{align*}
            \beta_k&= 2\sup_n \left\{ \sup_{A\in \ca_{n+k}^\infty} (\bbp(A|\ca_1^n)-\bbp(A)) \right\} \\
            &=\sup_n \left\{ \sup \sum_{i=1}^I \sum_{j=1}^J | \bbp(A_i\cap B_j) - \bbp(A_i) \bbp(B_j) |\right\}
        \end{align*}
        where the last supremum is over all finite $\ca_1^n$-measurable partitions $(A_1,\ldots,A_I)$ and all finite $\ca_{n+k}^\infty$-measurable partitions $(B_1,\ldots,B_J)$. 
        \item Let $\mathcal{G}$ and $\mathcal{H}$ be sub-$\sigma$-fields of $\mathcal{F}$. Then
        \[
            \alpha(\mathcal{G}, \mathcal{H}) := \sup_{G\in\mathcal{G}, H\in\mathcal{H}} \abs{\bbp(G \cap H) - \bbp(G) \bbp(H)}
        \]
        is called the \emph{strong mixing coefficient}, and we say that $(Z_t)_{t\in\bbz}$ is \emph{strongly mixing} if 
        \[
            \alpha_k := \sup_t \alpha(\ca^t_{-\infty}, \ca^\infty_{t+k}) \xrightarrow{} 0 
            \quad\text{as } k \to \infty.
        \]
    \end{enumerate}
\end{definition}

Absolute regularity implies strong mixing, but the other implication does not hold (see, for instance, \citealt{bradley80}). 
Sometimes, mixing conditions are too restrictive. Rather, it might be sufficient if the process at hand is only ``approximately mixing''.

\begin{definition}
\begin{enumerate}[label=(\alph*)]
    \item Let $(Z_t)_{t\in\bbz}$ be a real-valued stationary process. We call a sequence $(Y_t)_{t\in\bbz}$ a \emph{functional of} $(Z_t)_{t\in\bbz}$ if there is a measurable function $f$ defined on $\bbr^\bbz$ such that
    \begin{align*}
        Y_t=f((Z_{t+k})_{k\in\bbz}).
    \end{align*}
    \item Let $(Y_t)_{t\in\bbz}$ be a functional of $(Z_t)_{t\in\bbz}$ and let $r\geq 1$. Suppose that $(a_k)_{k\in\bbn_0}$ are constants with $a_k\to 0$.  We say that $(Y_t)_{t\in\bbz}$ satisfies the \emph{$r$-approximating condition} or that it is an $r$\emph{-approximating functional}, if 
    \begin{align} \label{approxatzero}
        \bbe\norm{Y_0-\bbe(Y_0|\ca_{-k}^k)}^r \leq a_k.
    \end{align}
    for all $k \in \bbn_0$.
    The sequence $(a_k)_{k \in \bbn_0}$ of approximating constants is said to be \emph{of size} $- \lambda$ if $a_k = \mathcal{O}(k^{-\lambda-\varepsilon})$ for some $\varepsilon > 0$. 
\end{enumerate}
\end{definition}

This kind of weak dependence is sometimes referred to by ``($L_p$-)near-epoch dependence'' in the literature. 
$r$-approximating functionals occur naturally in the theory of dynamical systems. Examples are given in \cite[Example~1.2--1.10]{borovkovaetal}, for instance.

We state the following definition in the multivariate version. \citet{schn_sil_mar_24} make use of the results of \citet{borovkovaetal}, where Sections~1--5 are all written down for the one-dimensional case. However, in Section~6 of \citet{borovkovaetal}, it is stated that all the previous results hold true in a multidimensional setting as well. The same applies to the results of \citet{schn_sil_mar_24} that we want to employ in the proof of Theorem~\ref{theorem: asymp distribution OP probs case 1}.

\begin{definition}[\protect{\citealp[Definition 2.10]{borovkovaetal}}]\label{definition: p-continuity}
For $m \geq 2$ and $q \geq 1$, let $h : \bbr^{m-1} \to \bbr^q$ be a measurable function and let $F$ be some distribution on $\bbr$. We say that $h$ is \emph{$p$-continuous with respect to $F$} if there exists a function $\phi : (0, \infty) \to (0, \infty)$ with $\phi(\varepsilon) = o(1)$ as $\varepsilon \to 0$ such that
\begin{equation}
    \bbe \left( \norm{h(Y) - h(Y^\prime)}^p \indicator{\norm{X - X^\prime} \leq \varepsilon}\right) \leq \phi(\varepsilon)
\end{equation}
holds for all random variables $Y$ and $Y^\prime$ with distribution $F$. If the choice of the underlying distribution $F$ is clearly understood, we simply say that $h$ is \emph{$p$-continuous}.
\end{definition}

The following theorem is a CLT for partial sums of functionals of absolutely regular processes derived by \citet{schn_sil_mar_24}. It is an extension of \citet[Theorem~4]{borovkovaetal} in terms of the posed summability condition.
For the readers' convenience, we state the multivariate version. 
Furthermore, note that since we are considering data in the respective statistics, for simplicity, we omit observations indexed by $t \leq 0$ in the following.

\begin{theorem}[\protect{\citealt[Theorem~
B.5]{schn_sil_mar_24}}]\label{theorem: Theorem 4}
    For $q \geq 1$, let $(\bfY_t)_{t\in\bbn}$, $\bfY_t = (Y_{t,1}, \dots, Y_{t,q})^\top$, be a 1-approximating functional with approximating constants $(a_k)_{k \in \bbn_0}$ of an absolutely regular process with mixing coefficients $(\beta_k)_{k\in\bbn_0}$. Furthermore, suppose that $\bbe \bfY_0 = \mathbf{0}$, all $(4+\delta)$-moments $\bbe \abs{Y_{0,i}}^{4+\delta} < \infty$ exist and
    \begin{equation}\label{eq: CLT Schnurr et al condition 1}
        \sum^{\infty}_{k=m} k^2 \brackets{a_k^{\frac{\delta}{3+\delta}} + \beta_k^{\frac{\delta}{4+\delta}}} < \infty
    \end{equation}
    for some $\delta > 0$ and a fixed integer $m \geq 0$. Then, as $n \to \infty$,
    \begin{equation*}
        \frac{1}{\sqrt{n}} \sum^n_{t=1} \bfY_t \xrightarrow{d} N(\mathbf{0}, \Sigma),
    \end{equation*}
    where $\Sigma=\bbe (\bfY_0 \bfY_0^\top) + 2 \sum^{\infty}_{k=1} \bbe (\bfY_0 \bfY_k^\top)$ converges absolutely. %In case $\sigma^2=0$, we adopt the convention that $N(0,0)$ denotes the point mass at the origin. 
    If $\bfY_0$ is bounded, the CLT continues to hold if \eqref{eq: CLT Schnurr et al condition 1} is replaced by the condition 
    \begin{equation}\label{eq: CLT Schnurr et al condition 2}
        \sum^{\infty}_{k=m} k^2 (a_k + \beta_k) < \infty.
    \end{equation}
\end{theorem}

%%%%%%%%%%%%%%%%%%%%
\section{Proofs}
\label{section: proofs}

%%%%%%%%%%% Proofs of Section 2
\subsection{Proof of Proposition~\ref{prop:p-continuous}}
\label{proof prop: p-continuous}
Let $\bfY$ and $\bfY^\prime$ denote two $(m-1)$-dimensional random vectors with the same distribution as $\bfDelta_0$. Every norm on $\bbr^{m!}$ is equivalent, so we use the taxi-norm $\norm{\cdot}_1$ here. Since $|\indicator{\Pi_{\bfDeltai}(\bfY)=\pi_i}-\indicator{\Pi_{\bfDeltai}(\bfY^\prime)=\pi_i}|\in \{0,1\}$ for all $i \in \{1, \dots, m!\}$, it follows
\begin{align}
    \bbe \left( \norm{h(\bfY) - h(\bfY^\prime)}^p \indicator{\norm{\bfY - \bfY^\prime}\leq \varepsilon}\right) 
    \leq m! \cdot \bbe\left(\indicator{\norm{\bfY - \bfY^\prime}\leq \varepsilon}\right) = m! \cdot \bbp(\norm{\bfY - \bfY^\prime}\leq \varepsilon) =: \phi(\varepsilon). \label{eq: phi specified}
\end{align}
The claim of Proposition~\ref{prop:p-continuous} follows, since $m!$ is a fixed factor and $\bbp(\norm{\bfY - \bfY^\prime}\leq \varepsilon) \to 0$ as $\varepsilon \to 0$.

\subsection{Proof of Theorem~\ref{theorem: asymp distribution OP probs case 1}}
\label{proof theorem: asymp distribution OP probs case 1}
The proof of \citet[Lemma~2.6]{schn_sil_mar_24} yields that $\bfDelta_t$ satisfies the 1-approximating condition for $k \geq m-2$ with approximating constants
\[
    2 \sum^{m-1}_{i=0} a_{k-i};
\]
also see our considerations before Proposition~\ref{prop:p-continuous}.
Since $|h|$ is bounded by 1, from \citet[Prop.\ 2.11]{borovkovaetal} it follows that $h(\bfDelta_t)$ satisfies the 1-approximating condition for $k \geq m-2$ with approximating constants 
\[
    a_k^\prime = \phi\left(\sqrt{4 \sum^{m-2}_{i=0} a_{k-i}}\right) + 2\cdot \sqrt{4 \sum^{m-2}_{i=0} a_{k-i}}.
\]
Hence, $h(\bfDelta_t)-\bfp$ is a bounded $1$-approximating functional and has mean zero. Note that the approximating constants as well as the absolute regularity coefficients are non-negative by definition, so with regard to the summability condition required by Theorem~\ref{theorem: Theorem 4}, it holds
\begin{align*}
    \sum^\infty_{k=m-2} k^2 \brackets{a_k^\prime + \beta_k} 
    %&= \sum^\infty_{k=m-2} k^2 \phi\left(2\sqrt{\sum^{m-2}_{i=0} a_{k-i}}\right) + 4 \sum^\infty_{k=m-2} k^2 \sqrt{\sum^{m-2}_{i=0} a_{k-i}} + \sum^\infty_{k=m-2} k^2 \beta_k \\
    &\leq \sum^\infty_{k=m-2} k^2 \phi\left(2\sqrt{\sum^{m-2}_{i=0} a_{k-i}}\right) + 4 \sum^\infty_{k=m-2} k^2 \sum^{m-2}_{i=0} \sqrt{a_{k-i}} + \sum^\infty_{k=m-2} k^2 \beta_k \\
    %&= \sum^\infty_{k=m-2} k^2 \phi\left(2\sqrt{\sum^{m-2}_{i=0} a_{k-i}}\right) + 4 \sum^{m-2}_{i=0} \sum^\infty_{k=m-2} k^2 \sqrt{a_{k-i}} + \sum^\infty_{k=m-2} k^2 \beta_k \\
    &= \sum^\infty_{k=m-2} k^2 \phi\left(2\sqrt{\sum^{m-2}_{i=0} a_{k-i}}\right) + 4 \sum^{m-2}_{i=0} \sum^\infty_{k=m-2-i} (k+i)^2 \sqrt{a_{k}} + \sum^\infty_{k=m-2} k^2 \beta_k,
\end{align*}
where we have applied the elementary inequality $\sqrt{a+b} \leq \sqrt{a} + \sqrt{b}$ for $a,b \geq 0$ (see \citealt[p.~157]{loe_77}).
Since
\begin{align*}
    \sum^\infty_{k=m-2-i} (k+i)^2 \sqrt{a_{k}} 
    &= \sum^\infty_{k=m-2-i} k^2 \sqrt{a_k} + 2i \sum^\infty_{k=m-2-i} k \sqrt{a_k} + i^2 \sum^\infty_{k=m-2-i} \sqrt{a_k} \\
    &\leq (1+2i+i^2) \sum^\infty_{k=0} k^2 \sqrt{a_k} < \infty
\end{align*}
by assumption, we see that 
\[
    \sum^\infty_{k=m-2} k^2 \brackets{a_k^\prime + \beta_k} < \infty.
\]
Thus, from Theorem~\ref{theorem: Theorem 4} it follows that 
\[
    \frac{1}{\sqrt{n}} \sum^n_{t=1} (h(\bfDelta_t) - \bfp) \xrightarrow{d} N(0, \Sigma)
\]
as $n \to \infty$, where $\Sigma$ is given by \eqref{eq: long-run covariance matrix}. This completes the proof of Theorem~\ref{theorem: asymp distribution OP probs case 1}.

\subsection{Remark}
\label{appendix remark}
In Theorem \ref{theorem: asymp distribution OP probs case 1}, we did not pose the conditions on the observed process itself, but on the process consisting of the increments. However, posing the approximating condition and stationarity on the process $(X_t)_{t\in\bbn_0}$ itself, it follows
\begin{align*}
    \bbe|\Delta_0 - \bbe(\Delta_0|\ca^k_{-k})|
    &\leq \bbe|X_0 - \bbe(X_0|\ca^k_{-k})| + \bbe|X_{-1} - \bbe(X_{-1}|\ca^k_{-k})| \\
    &\leq \bbe|X_0 - \bbe(X_0|\ca^k_{-k})| + 2 \cdot \bbe|X_{-1} - \bbe(X_{-1}|\ca^{k-2}_{-k})|,
\end{align*}
where we have used  
\citet[Theorem~10.28]{davidson}. Since 
$\bbe|X_{-1} - \bbe(X_{-1}|\ca^{k-2}_{-k})| = \bbe|X_0 - \bbe(X_0|\ca^{k-1}_{-(k-1)})|$ by \citet[Lemma~2.5]{schn_sil_mar_24}, application of the $1$-approximation condition yields
\[
    \bbe|\Delta_0 - \bbe(\Delta_0|\ca^k_{-k})| \leq a_k + 2 a_{k-1},
\]
where $(a_k)_{k\geq0}$ denote the approximating constants of $(X_t)_{t \in\bbn_0}$. Hence, $(\Delta_t)_{t\in\bbn_0}$ satisfies the $1$-approximating condition for $k \geq 1$. Since stationarity implies increment stationarity, we obtain that $(\Delta_t)_{t\in\bbn_0}$ is also stationary. Continuing in the same fashion as in the proof of Theorem~\ref{theorem: asymp distribution OP probs case 1}, for $k \geq m$, we obtain 
\[
    a_k^\prime := \phi\brackets{\sqrt{4 \sum^{m-1}_{i=0} (a_{k-i} + 2 a_{k-1-i})}} + 2 \cdot \sqrt{4 \brackets{\sum^{m-1}_{i=0} (a_{k-i} + 2 a_{k-1-i})}}
\]
as the approximating constants of $h(\overline{\Delta_t}) - \bfp$,
which implies that \eqref{eq: asymp distribution OP probs case 1 phi} has to be changed to 
\begin{equation}
\label{eq: asymp distribution OP probs case 1 phi altered}
    \sum^\infty_{k=m} k^2 \phi\brackets{2 \cdot\sqrt{\sum^{m-1}_{i=0} (a_{k-i} + 2 a_{k-1-i})}} < \infty
\end{equation}
for Theorem~\ref{theorem: Theorem 4} to be still applicable. However, since \eqref{eq: asymp distribution OP probs case 1 phi altered} implies \eqref{eq: asymp distribution OP probs case 1 phi}, 
this shows that posing the conditions on the increment process can indeed be understood as the weaker requirement for the claim to hold.

%%%%%%%%%%%%%%%%%%%%% Proofs of Section 3
\subsection{Proof of Proposition~\ref{prop: cov MA(q) m2}}
\label{proof prop: cov MA(q) m2}
First, let us consider the diagonal entries. For $m=2$, there are only two OPs, namely the increasing OP $(1,2)$ and the decreasing OP $(2,1)$, both with probability $\tfrac{1}{2}$ according to \citet[Proposition~7]{bandt07}. Starting with the increasing OP, we have
\begin{equation*}
\label{eq: MAq m2 step1}
    \sigma_{11} = \frac{1}{4} 
    + \sum^{q+1}_{k=1} \Bigl(2 \cdot \bbp(X_0 < X_1, X_k < X_{k+1}) - \frac{1}{2}\Bigr).
\end{equation*}
By definition, $X_0 = \sum^{q-1}_{l=0} \varepsilon_{-l} + \varepsilon_{-q}$ and $X_1 = \varepsilon_1 + \sum^q_{l=1} \varepsilon_{1-l} = \varepsilon_1 + \sum^{q-1}_{l=0} \varepsilon_{-l}$, so $X_0 < X_1$ is equivalent to $\varepsilon_{-q} < \varepsilon_1$; analogously, $X_k < X_{k+1}$ is equivalent to $\varepsilon_{k-q} < \varepsilon_{k+1}$ (cf. \citealp{bandt07}). Hence, 
\begin{align}
\label{eq: MAq m2 step11}
    \sigma_{11} = \frac{1}{4} + \sum^{q+1}_{k=1} \Bigl(2 \cdot \bbp(\varepsilon_{-q} < \varepsilon_1, \varepsilon_{k-q} < \varepsilon_{k+1}) - \frac{1}{2}\Bigr). 
\end{align}
Note that the inequalities inside the probability refer to mutually independent noise terms for $1 \leq k \leq q$, so by the i.i.d.\ assumption, it follows that
\begin{equation}
\label{eq: MAq m2 step2}
    \bbp(\varepsilon_{-q} < \varepsilon_1, \varepsilon_{k-q} < \varepsilon_{k+1}) = \bbp(\varepsilon_{-q} < \varepsilon_1) \bbp(\varepsilon_{k-q} < \varepsilon_{k+1}) = \bbp(\varepsilon_1<\varepsilon_2)^2 = \frac{1}{4}, \quad \text{for } 1 \leq k \leq q.
\end{equation}
Hence, the summands for $1 \leq k \leq q$ in \eqref{eq: MAq m2 step2} vanish.
Furthermore, the case $k=q+1$ yields
\begin{equation}
\label{eq: MAq m2 step3}
    \bbp(\varepsilon_{-q} < \varepsilon_1 < \varepsilon_{q+2}) = \bbp(\varepsilon_1 < \varepsilon_2 < \varepsilon_3) = \frac{1}{6}.
\end{equation}
Combining \eqref{eq: MAq m2 step2} and \eqref{eq: MAq m2 step3}, we conclude
\[
    \sigma_{11} = \frac{1}{4} %+ \sum^{q}_{k=1} \Bigl(2 \cdot \frac{1}{4} - \frac{1}{2}\Bigr) 
    + \Bigl(2 \cdot \frac{1}{6} - \frac{1}{2}\Bigr) = \frac{1}{12}.
\]
Similarly, one can show $\sigma_{22} = \tfrac{1}{12}$. It remains to derive the off-diagonal entry $\sigma_{12} = \sigma_{21}$, where we proceed analogously:
\begin{align*}
    \sigma_{12} &= - \frac{1}{4} + \sum^{q+1}_{k=1} \Bigl(\bbp(X_0<X_1, X_k>X_{k+1}) + \bbp(X_0>X_1, X_k<X_{k+1}) - \frac{1}{2}\Bigr) \\
    &= - \frac{1}{4} + \sum^{q+1}_{k=1} \Bigl(\bbp(\varepsilon_{-q}<\varepsilon_1, \varepsilon_{k-q}>\varepsilon_{k+1}) + \bbp(\varepsilon_{-q}>\varepsilon_1, \varepsilon_{k-q}<\varepsilon_{k+1}) - \frac{1}{2}\Bigr) \\
    &= - \frac{1}{4} + \sum^{q}_{k=1} \Bigl(\bbp(\varepsilon_{-q}<\varepsilon_1) \bbp(\varepsilon_{k-q}>\varepsilon_{k+1}) + \bbp(\varepsilon_{-q}>\varepsilon_1) \bbp(\varepsilon_{k-q}<\varepsilon_{k+1}) - \frac{1}{2}\Bigr) \\
    &\hspace{20mm} + \Bigl(\bbp(\varepsilon_{-q}<\varepsilon_1, \varepsilon_{1}>\varepsilon_{q+2}) + \bbp(\varepsilon_{-q}>\varepsilon_1, \varepsilon_{1}<\varepsilon_{q+2}) - \frac{1}{2}\Bigr)\\
    &= - \frac{1}{4} + \sum^{q}_{k=1} \Bigl(\frac{1}{4} + \frac{1}{4} - \frac{1}{2}\Bigr) + \Bigl(\bbp(\varepsilon_1 < \varepsilon_2, \varepsilon_2 > \varepsilon_3) + \bbp(\varepsilon_1 > \varepsilon_2, \varepsilon_2 < \varepsilon_3) - \frac{1}{2} \Bigr).
\end{align*}
The events in the remaining two probabilities correspond to the events that $(\varepsilon_1, \varepsilon_2, \varepsilon_3)^\top$ has OP $(1, 3, 2)$ or $(2, 3, 1)$, and that $(\varepsilon_1, \varepsilon_2, \varepsilon_3)^\top$ has OP $(2, 1, 3)$ or $(3, 1, 2)$, respectively. Due to the i.i.d.\ assumption, the OP distributions are uniform, so each of these probabilities equals~$\tfrac{1}{3}$, which yields $\sigma_{12} = - \tfrac{1}{12}$. Altogether the assertion of Proposition~\ref{prop: cov MA(q) m2} follows.

%%%%%%%%%%%%%%%%%%%%%%
\subsection{Proof of Proposition~\ref{prop: Sigma_ct}}
\label{proof prop: Sigma_ct}

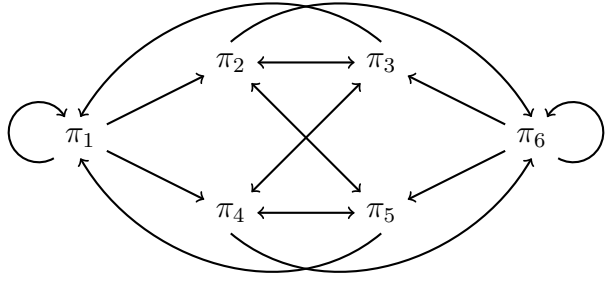
\begin{figure}[t]
    \centering
    \begin{tikzpicture}[squarednode/.style={rectangle, %draw=black, 
    thick, minimum size=5mm}]
        %Nodes
        \node[squarednode] (pi1) at (0,0) {$\pi_1$};
        \node[squarednode] (pi2) at (2, 1) {$\pi_2$};
        \node[squarednode] (pi3) at (4, 1) {$\pi_3$};
        \node[squarednode] (pi4) at (2, -1) {$\pi_4$};
        \node[squarednode] (pi5) at (4, -1) {$\pi_5$};
        \node[squarednode] (pi6) at (6, 0) {$\pi_6$};
        %Lines
        \draw[->, thick] (pi1.south west) arc[start angle=300, end angle=30, radius=0.4cm];
        \draw[->, thick] (pi1) -- (pi2);
        \draw[->, thick] (pi1) -- (pi4);
        
        \draw[<->, thick] (pi2) -- (pi3);
        \draw[<->, thick] (pi2) -- (pi5);
        \draw[<->, thick] (pi3) -- (pi4);
        \draw[<->, thick] (pi4) -- (pi5);

        \draw [->, thick] (pi2.north) to [out=40,in=120] (pi6.north);
        \draw [->, thick] (pi4.south) to [out=-40,in=-120] (pi6.south);
        \draw [->, thick] (pi3.north) to [out=140,in=60] (pi1.north);
        \draw [->, thick] (pi5.south) to [out=-140,in=-60] (pi1.south);

        \draw[->, thick] (pi6.south east) arc[start angle=-120, end angle=150, radius=0.4cm];
        \draw[->, thick] (pi6) -- (pi3);
        \draw[->, thick] (pi6) -- (pi5);
    \end{tikzpicture}
    \caption{Allowed transitions for OPs of length $m=3$.}
    \label{fig: OP transitions}
\end{figure}

Recall that $q=1-p$. Furthermore, recall that we denote the events $X_i < X_j$ and $X_i > X_j$ by $C_{ij}^+$ and $C_{ij}^-$, respectively, i.e., the superscripts $+,-$ indicate whether the order is increasing or decreasing. In order to simplify notation, we omit the subscript whenever the information obtained from it is either not necessary or given by the context. The probabilities of OPs of length~3, which we need to compute the terms of the form $p_i (\delta_{ij} - p_j)$ for $1 \leq i,j \leq 6$, are given in Table~\ref{tab: OP probs coin toss m3} together with the coin tosses that lead to the respective probabilities.
For $k \geq 2$, 
\[
    \bbp(\Pi_0 = \pi_i, \Pi_k=\pi_j) = \bbp(\Pi_0 = \pi_i) \bbp(\Pi_k=\pi_j) = p_i \cdot p_j,
\]
since the tosses needed to determine the OP $\Pi_k$ are independent of all the tosses involving $X_0$ and $X_1$. Hence, it remains to determine the probabilities $\bbp(\Pi_0 = \pi_i, \Pi_1=\pi_j)$. Note that $\Pi_0$ and $\Pi_1$ share the components $X_1$ and $X_2$. As a consequence, not every OP can follow every other OP. Figure~\ref{fig: OP transitions} shows all allowed transitions. In order to compute the probabilities, we must consider all OPs of length~4 such that the first and last three components satisfy the respective OPs of length~3. The latter probabilities are stated in Table~\ref{tab: OP probs coin toss m4}. 

\begin{table}[t]
    \label{tab: OP probs coin toss m3}
    \caption{OPs of length $m=3$, the required coin tosses and the resulting probabilities depending on the parameter $p$ and $q=1-p$.}
    \centering
    \makegapedcells
    \begin{tabular}{cclc}
        \toprule
        Index & OP & Tosses & Probability \\ \midrule
        1 & (1, 2, 3) & $C_{12}^+, C_{23}^+$ & $p^2$ \\ 
        2 & (1, 3, 2) & $C_{12}^+, C_{23}^-, C_{13}^+$  & $p^2 q$ \\
        3 & (2, 1, 3) & $C_{12}^-, C_{23}^+, C_{13}^+$ & $p^2 q$ \\
        4 & (2, 3, 1) & $C_{12}^+, C_{23}^-, C_{13}^-$ & $p q^2$ \\
        5 & (3, 1, 2) & $C_{12}^-, C_{23}^+, C_{13}^-$ & $p q^2$ \\
        6 & (3, 2, 1) & $C_{12}^-, C_{23}^-$ & $q^2$ \\
        \bottomrule
    \end{tabular}
\end{table}

\begin{table}[t]
    \label{tab: OP probs coin toss m4}
    \caption{OPs of length $m=4$, the required coin tosses and the resulting probabilities depending on the parameter $p$ and $q=1-p$.}
    \centering
    \makegapedcells
    \begin{tabular}{cclc}
        \toprule
        Index & OP & Tosses & Probability  \\ \midrule
        1 & (1, 2, 3, 4) & $C_{12}^+, C_{23}^+,\ C_{34}^+$ & $p^3$ \\
        2 & (1, 2, 4, 3) & $C_{12}^+, C_{23}^+,\ C_{34}^-, C_{24}^+$ & $p^3 \cdot q$ \\
        3 & (1, 3, 2, 4) & $C_{12}^+, C_{23}^-, C_{13}^+,\ C_{34}^+, C_{24}^+$ & $p^4 \cdot q$ \\
        4 & (1, 3, 4, 2) & $C_{12}^+, C_{23}^+,\ C_{34}^-, C_{24}^-, C_{14}^+$ & $p^3 \cdot q^2$ \\
        5 & (1, 4, 2, 3) & $C_{12}^+, C_{23}^-, C_{13}^+,\ C_{34}^+, C_{24}^-$ & $p^3 \cdot q^2$ \\
        6 & (1, 4, 3, 2) & $C_{12}^+, C_{23}^-, C_{13}^+,\ C_{34}^-, C_{14}^+$ & $p^3 \cdot q^2$ \\
        7 & (2, 1, 3, 4) & $C_{12}^-, C_{23}^+, C_{13}^+,\ C_{34}^+$ & $p^3 \cdot q$ \\
        8 & (2, 1, 4, 3) & $C_{12}^-, C_{23}^+, C_{13}^+,\ C_{34}^-, C_{24}^+, C_{14}^+$ & $p^4 \cdot q^2$ \\
        9 & (2, 3, 1, 4) & $C_{12}^+, C_{23}^-, C_{13}^-,\ C_{34}^+, C_{23}^+$ & $p^3 \cdot q^2$ \\
        10 & (2, 3, 4, 1) & $C_{12}^+, C_{23}^+,\ C_{34}^-, C_{24}^-, C_{14}^-$ & $p^2 \cdot q^3$ \\
        11 & (2, 4, 1, 3) & $C_{12}^+, C_{23}^-, C_{13}^-,\ C_{34}^+, C_{24}^-, C_{14}^+$ & $p^3 \cdot q^3$ \\
        12 & (2, 4, 3, 1) & $C_{12}^+, C_{23}^-, C_{13}^+,\ C_{34}^-, C_{14}^-$ & $p^2 \cdot q^3$ \\
        13 & (3, 1, 2, 4) & $C_{12}^-, C_{23}^+, C_{13}^-,\ C_{34}^+, C_{14}^+$ & $p^3 \cdot q^2$ \\ 
        14 & (3, 1, 4, 2) & $C_{12}^-, C_{23}^+, C_{13}^+,\ C_{34}^-, C_{24}^+, C_{14}^-$ & $p^3 \cdot q^3$ \\
        15 & (3, 2, 1, 4) & $C_{12}^-, C_{23}^-,\ C_{34}^+, C_{24}^+, C_{14}^+$ & $p^3 \cdot q^2$ \\
        16 & (3, 2, 4, 1) & $C_{12}^-, C_{23}^+, C_{13}^+,\ C_{34}^-, C_{24}^-$ & $p^2 \cdot q^3$ \\
        17 & (3, 4, 1, 2) & $C_{12}^+, C_{23}^-, C_{13}^-,\ C_{34}^+, C_{24}^-, C_{14}^-$ & $p^2 \cdot q^4$ \\
        18 & (3, 4, 2, 1) & $C_{12}^+, C_{23}^-, C_{13}^-,\ C_{34}^-$ & $p \cdot q^3$ \\
        19 & (4, 1, 2, 3) & $C_{12}^-, C_{23}^+, C_{13}^-,\ C_{34}^+, C_{14}^-$ & $p^2 \cdot q^3$ \\
        20 & (4, 1, 3, 2) & $C_{12}^-, C_{23}^+, C_{13}^-,\ C_{34}^-, C_{24}^+$ & $p^2 \cdot q^3$ \\
        21 & (4, 2, 1, 3) & $C_{12}^-, C_{23}^-,\ C_{34}^+, C_{24}^+, C_{14}^-$ & $p^2 \cdot q^3$ \\
        22 & (4, 2, 3, 1) & $C_{12}^-, C_{23}^+, C_{13}^-,\ C_{34}^-, C_{24}^-$ & $p \cdot q^4$ \\
        23 & (4, 3, 1, 2) & $C_{12}^-, C_{23}^-,\ C_{34}^+, C_{24}^-$ & $p \cdot q^3$ \\
        24 & (4, 3, 2, 1) & $C_{12}^-, C_{23}^-,\ C_{34}^-$ & $q^3$ \\
        \bottomrule
    \end{tabular}
\end{table}

Let us explain the procedure with a few examples. It holds that
\[
    \bbp(\Pi_0 =\pi_1, \Pi_1 = \pi_1) = \bbp(X_0 < X_1 < X_2 < X_3) = p^3,
\]
because we need three coin tosses $C^+$ and the remaining orderings are determined by transitivity. Next, we consider the transition from $\pi_1$ to $\pi_2$. The only OP of length 4 such that $\Pi(\bfX_0) = \pi_1$ and $\Pi(\bfX_1) = \pi_2$ are satisfied is $(1, 2, 4, 3)$. It has probability $p^3 \cdot q$ because we need an additional coin toss $C_{24}^+$ to determine the order between $X_2$ and $X_4$. Finally, both OPs $(1, 3, 4, 2)$  and $(2, 3, 4, 1)$ result in the transition of $\pi_1$ to $\pi_4$, with probabilities $p^3\cdot q^2$ and $p^2 \cdot q^3$, respectively. Hence, $\bbp(\Pi_0 =\pi_1, \Pi_1 = \pi_4) = p^2 \cdot q^2$. Continuing in this fashion, we obtain the probability matrix
\begin{align*}
    \mathbf{P}_3(1) &:= \bigl(\bbp(\Pi_0 =\pi_i, \Pi_1 = \pi_j)\bigr)_{1 \leq i,j \leq 6} 
    =   \scalemath{0.85}{\begin{pmatrix}
        p^3 & p^3q & 0 & p^2q^2 & 0 & 0 \\
        0 & 0 & p^4q & 0 & p^3q^2& p^2q^2 \\
        p^3q & p^3q^2 & 0 & p^2q^3 & 0 & 0 \\
        0 & 0 & p^3q^2 & 0 & p^2q^3 & pq^3 \\
        p^2q^2 & p^2q^3 & 0 & pq^4 & 0 & 0 \\
        0 & 0 & p^2q^2 & 0 & pq^3 & q^3
    \end{pmatrix}},
\end{align*}
which leads to
\begin{align*}
    \Sigma &= \scalemath{0.85}{\begin{pmatrix}
        p^2(1-p^2) & -p^4q & -p^4q & -p^3q^2 & -p^3q^2 & -p^2q^2 \\
        -p^4q & p^2q(1-p^2q) & -p^4q^2 & -p^3q^3 & -p^3q^3 & -p^2q^3 \\
        -p^4q & -p^4q^2 & p^2q(1-p^2q) & -p^3q^3 & -p^3q^3 & -p^2q^3 \\
        -p^3q^2 & -p^3q^3 & -p^3q^3 & pq^2(1-pq^2) & -p^2q^4 & -pq^4 \\
        -p^3q^2 & -p^3q^3 & -p^3q^3 & -p^2q^4 & pq^2(1-pq^2) & -pq^4 \\
        -p^2q^2 & -p^2q^3 & -p^2q^3 & -pq^4 & -pq^4 & q^2(1-q^2)
    \end{pmatrix}} \\
    &\hspace{20mm} + \mathbf{P}_3(1) + \mathbf{P}_3(1)^\top - 2 \scalemath{0.85}{\begin{pmatrix}
        p^4 & p^4q & p^4q & p^3q^2 & p^3q^2 & p^2q^2 \\
        p^4q & p^4q^2 & p^4q^2 & p^3q^3 & p^3q^3 & p^2q^3 \\
        p^4q & p^4q^2 & p^4q^2 & p^3q^3 & p^3q^3 & p^2q^3 \\
        p^3q^2 & p^3q^3 & p^3q^3 & p^2q^4 & p^2q^4 & pq^4 \\
        p^3q^2 & p^3q^3 & p^3q^3 & p^2q^4 & p^2q^4 & pq^4 \\
        p^2q^2 & p^2q^3 & p^2q^3 & p^1q^4 & p^1q^4 & q^4
    \end{pmatrix}}. %\\
    % &= \begin{pmatrix}
    %     p^2(1+2p-3p^2) &p^3q(1-3p) &p^3q(1-3p) &p^2q^2(1-3p) &p^2q^2(1-3p) &-3p^2q^2 \\
    %     p^3q(1-3p) &p^2q(1-3p^2q) &p^3q(1-3pq) &-3p^3q^3 &p^2q^2(1-3pq) &p^2q^2(1-3q) \\
    %     p^3q(1-3p) &p^3q(1-3pq) &p^2q(1-3p^2q) &p^2q^2(1-3pq) &-3p^3q^3 &p^2q^2(1-3q) \\
    %     p^2q^2(1-3p) &-3p^3q^3 &p^2q^2(1-3pq) &pq^2(1-3pq^2) &pq^3(1-3pq) &pq^3(3p-2) \\
    %     p^2q^2(1-3p) &p^2q^2(1-3pq) &-3p^3q^3 &pq^3(1-3pq) &pq^2(1-3pq^2) &pq^3(3p-2) \\
    %     -3p^2q^2 &p^2q^2(1-3q) &p^2q^2(1-3q) &pq^3(3p-2) &pq^3(3p-2) &pq^2(4-3p)
    % \end{pmatrix}
\end{align*}
This expression further simplifies to the matrix given in \eqref{Sigma_ct}.

%%%%%%%%%%%%%%% Proofs of Section 4
\subsection{Proof of Lemma~\ref{lemmaH}}
\label{proof lemmaH}
The Delta method together with the first-order Taylor expansion \eqref{PE_Taylor1} immediately implies the following asymptotic normality of the permutation entropy:
\begin{equation*}
\sqrt{n} (H(\hat{\bfp}) - H(\bfp)) \xrightarrow{d} N(0, \sigma^2_{PE}),
\quad\text{where}\quad
\sigma^2_{PE} = \sum_{i,j=1}^{m!} (1+ \log{p_i})(1+\log{p_j})\, \sigma_{ij}.
\end{equation*}
Analogously, evaluating $H\big((\hat{\bfp}+\bfu)/2\big)$ around $(\bfp+\bfu)/2$, we get
$$
H\big((\hat{\bfp}+\bfu)/2\big)\ \approx\ H\big((\bfp+\bfu)/2\big)
 - \sum_{k=1}^{m!} \frac{1+\log{\frac{p_{k}+u_k}{2}}}{2}\,\big(\hat{p}_k-p_{k}\big).
$$
Hence, the joint asymptotic distribution of~$H(\hat{\bfp})$ and~$H\big((\hat{\bfp}+\bfu)/2\big)$ is given by
\begin{equation*}
\sqrt{n}\, \left(\begin{array}{c} H(\hat{\bfp}) - H(\bfp)\\ H\big((\hat{\bfp}+\bfu)/2\big) - H\big((\bfp+\bfu)/2\big) \end{array}\right) \xrightarrow{d} N(\mathbf{0}, \Sigma^{(1)})
\text{ with }
\Sigma^{(1)}=B_1\,\Sigma\, B_1^\top,
\end{equation*}
where~$B_1$ abbreviates the Jacobian matrix
$$
B_1\ =\ -\left(\begin{array}{ccc}
1+\log{p_{1}} & \cdots & 1+\log{p_{m!}} \\
\frac{1}{2}\,\big(1+\log{\frac{p_{1}+1/m!}{2}}\big) & \cdots & \frac{1}{2}\,\big(1+\log{\frac{p_{m!}+1/m!}{2}}\big) \\
\end{array}\right).
$$
Evaluating the matrix product $B_1\,\Sigma\, B_1^\top$, we obtain 
\begin{eqnarray*}
\sigma^{(1)}_{11} &=& \textstyle
\sum_{i,j=1}^{m!} (1+ \log{p_i})(1+\log{p_j})\, \sigma_{ij},
\\
\sigma^{(1)}_{12} &=& \textstyle
\frac{1}{2}\,\sum_{i,j=1}^{m!} (1+ \log{p_i})(1+\log{\frac{p_j+1/m!}{2}})\, \sigma_{ij},
\\
\sigma^{(1)}_{22} &=& \textstyle
\frac{1}{4}\,\sum_{i,j=1}^{m!} (1+ \log{\frac{p_i+1/m!}{2}})(1+\log{\frac{p_j+1/m!}{2}})\, \sigma_{ij}.
\end{eqnarray*}
Recall that the $\sigma_{ij}$ from \eqref{eq: long-run covariance matrix} are given by
$$
\sigma_{ij} = p_i (\delta_{ij} - p_j) + \sum^\infty_{k=1} \Bigl(p_{ij}(k) + p_{ji}(k) - 2 p_i p_j\Bigr).
$$
Hence, using that $\sum_{j=1}^{m!} p_{ij}(k) = p_i = \sum_{j=1}^{m!} p_{ji}(k)$, we get
$$
\textstyle
\sum_{j=1}^{m!} \sigma_{ij}\ =\ p_i \cdot (1-\sum_{j=1}^{m!} p_j) + \sum^\infty_{k=1} \Bigl(\sum_{j=1}^{m!} p_{ij}(k) + \sum_{j=1}^{m!} p_{ji}(k) - 2 p_i \sum_{j=1}^{m!} p_j\Bigr)\ =\ 0.
$$
Furthermore, let $\bfp'$ abbreviate an arbitrary $m!$-dimensional pmf~vector. Then,
$$
\begin{array}{@{}rl}
\sum_{i,j=1}^{m!} \log{p_j'}\, \sigma_{ij}
\ =& \sum_{i=1}^{m!} \, p_i \big(\log{p_i'} - \sum_{j=1}^{m!} p_j \log{p_j'}\big)
\\[1ex]
&\ + \sum^\infty_{k=1} \Bigl(\sum_{j=1}^{m!} \log{p_j'}\cdot p_j + \sum_{j=1}^{m!} \log{p_j'}\, \cdot p_j - 2 \sum_{j=1}^{m!} \log{p_j'}\, p_j\Bigr)\ =\ 0.
\end{array}
$$
Hence, the entries of~$\Sigma^{(1)}$ can be further simplified, leading to the closed-form expressions in Lemma~\ref{lemmaH}, so the proof is complete.

\subsection{Proof of Lemma~\ref{lemma: det Sigma1 m2}}
\label{proof lemma: det Sigma1 m2}
To prove Lemma~\ref{lemma: det Sigma1 m2}, it suffices to show $\sigma_{ij} \cdot \sigma_{kl} = \sigma_{ik}\cdot \sigma_{jl}$ for $1 \leq i,j,k,l \leq 2$, where $\sigma_{ij}$ is defined by \eqref{eq: long-run covariance matrix}. Since $\sigma_{ij} = \sigma_{ji}$, it only remains to show that $\sigma_{11} \cdot \sigma_{22} =\sigma_{12}^2$.
    For future reference, let us state the following evident facts:
    \begin{align}
        1 &= p_1+p_2 \label{eq: det Sigma1 m2 1} \\
        \text{ and } \quad p_i &= p_{i1}(k) + p_{i2}(k) = p_{1i}(k) + p_{2i}(k) \quad \text{ for } i=1,2 \text{ and all } k \geq 1. \label{eq: det Sigma1 m2 2}
    \end{align}
    Absolute convergence of the respective sums (cf. Theorem~\ref{theorem: asymp distribution OP probs case 1}) and \eqref{eq: det Sigma1 m2 1} yields
    \begin{align*}
        \sigma_{11} \cdot \sigma_{22} 
        &= p_1(1-p_1)p_2(1-p_2) \\
        &\hspace{15mm}\textstyle + p_1(1-p_1) \sum_{k=1}^\infty 2(p_{22}(k)-p_2^2) + p_2(1-p_2) \sum_{k=1}^\infty 2(p_{11}(k)-p_1^2) \\
        &\hspace{15mm}\textstyle + \sum_{k,l=1}^\infty (2p_{11}(k)-2p_1^2)\cdot(2p_{22}(l)-2p_2^2) \\
        &\textstyle
        = (p_1p_2)^2 + 2p_1p_2 \sum_{k=1}^\infty \underbrace{(p_{11}(k)+p_{22}(k)-p_1^2-p_2^2)}_{\text{(I)}} \\
        &\hspace{15mm}\textstyle + \sum_{k,l=1}^\infty \underbrace{(2p_{11}(k)-2p_1^2)}_{\text{(II)}}\cdot(2p_{22}(l)-2p_2^2).
    \end{align*}
    First, we consider (I). By repeated application of \eqref{eq: det Sigma1 m2 1} and \eqref{eq: det Sigma1 m2 2}, we obtain
    \begin{align*}
        \text{(I)} &= (p_1-p_{12}(k))+(p_2-p_{21}(k)) - p_1(1-p_2) - p_2(1-p_1) \\
        &= 1-p_{12}(k)-p_{21}(k) + 2p_1p_2-1 \\
        &=-(p_{12}(k)+p_{21}(k)-2p_1p_2).
    \end{align*}
    Similarly, it follows that
    \begin{align*}
        \text{(II)} &= p_{11}(k)+p_{11}(k)-2p_1(1-p_2) \\
        &= (p_1-p_{12}(k)) + (p_1 - p_{21}(k)) -2p_1 + 2p_1p_2 \\
        &= -(p_{12}(k) + p_{21}(k) - 2p_1p_2),
    \end{align*}
    and analogously $2p_{22}(l)-2p_2^2 = -(p_{12}(l) + p_{21}(l) - 2p_1p_2)$. Combining these, we conclude
    \begin{align*}
        \sigma_{11}\sigma_{22} 
        &\textstyle
        = (-p_1p_2)^2 - 2p_1p_2 \sum_{k=1}^\infty (p_{12}(k)+p_{21}(k)-2p_1p_2) \\
        &\hspace{15mm} \textstyle + \sum_{k,l=1}^\infty (-1)^2 (p_{12}(k)+p_{21}(k)-2p_1p_2)(p_{12}(l)+p_{21}(l)-2p_1p_2) \\
        &= \sigma_{12}^2,
    \end{align*}
    so the proof of Lemma~\ref{lemma: det Sigma1 m2} is complete.

\subsection{Proof of Theorem~\ref{theorem: asymp distr H D}}
\label{proof theorem: asymp distr H D}
We define the function $g(x,y)= \big(x,\ y - \tfrac{1}{2}\,x - \tfrac{1}{2}\,\log{m!}\big)^\top$ such that $\big(H(\bfp), D(\bfp)\big)^\top$ equals $g\big(H(\bfp), H((\bfp+\bfu)/2)\big)$. 
The function $g(x,y)$ has the Jacobian $J_g(x,y) = \big(\begin{smallmatrix} 1 & 0 \\ - 1/2 & 1 \end{smallmatrix}\big) =: B_2$, so applying the Delta method to \eqref{HH_joint_asymptotics} in Lemma~\ref{lemmaH}, it follows that
\begin{eqnarray*}
%\label{HD_joint_asymptotics}
\sqrt{n}\, \left(\begin{array}{c} H(\hat{\bfp}) - H(\bfp)\\ D(\hat{\bfp}) - D(\bfp) \end{array}\right) &\xrightarrow{d}& N(\mathbf{0}, \Sigma^{(2)})
\text{ with }
\\
\nonumber
\Sigma^{(2)}=B_2\,\Sigma^{(1)}\, B_2^\top
&=&
\left(\begin{array}{cc}
\sigma^{(1)}_{11} & \sigma^{(1)}_{12}-\sigma^{(1)}_{11}/2 \\
\sigma^{(1)}_{12}-\sigma^{(1)}_{11}/2 & \sigma^{(1)}_{22}-\sigma^{(1)}_{12}+\sigma^{(1)}_{11}/4 \\
\end{array}\right).
\end{eqnarray*}
The resulting explicit expressions are those stated in Theorem~\ref{theorem: asymp distr H D}, which follows by using that $\log{\frac{p_j+1/m!}{2}}-\log{p_j} = \log{\frac{p_j+1/m!}{2p_j}}$. So the proof of Theorem~\ref{theorem: asymp distr H D} is complete.

\subsection{Proof of Theorem~\ref{theorem: asymp distr H C}}
\label{proof theorem: asymp distr H C}
We define the function $h(x,y) = \big(H_0\, x,\ H_0 D_0\, x\,(y - \tfrac{1}{2}\,x - \tfrac{1}{2}\,\log{m!})\big)^\top$, which maps $\big(H(\bfp),\, H((\bfp+\bfu)/2)\big)$ onto $\big(H_0\,H(\bfp),\, C(\bfp)\big)$. As the corresponding Jacobian, we get
$$
B_3 := J_h\big(H(\bfp),\, H((\bfp+\bfu)/2)\big)
%\ =\ 
%\left(\begin{array}{cc}
%H_0 & 0 \\
%%H_0 D_0\, (H((\bfp+\bfu)/2)-H(\bfp) - \tfrac{1}{2}\,\log{d!}) & H_0 D_0\, H(\bfp) \\
%%H_0 D_0\, (H((\bfp+\bfu)/2)-\tfrac{1}{2}\,H(\bfp) - \tfrac{1}{2}\,\log{d!} - \tfrac{1}{2}\,H(\bfp)) & H_0 D_0\, H(\bfp) \\
%H_0 D_0\, (D(\bfp) - \tfrac{1}{2}\,H(\bfp)) & H_0 D_0\, H(\bfp) \\
%\end{array}\right).
\ =\ 
H_0\left(\begin{array}{cc}
1 & 0 \\
D_0\, (D(\bfp) - \tfrac{1}{2}\,H(\bfp)) & D_0\, H(\bfp) \\
\end{array}\right).
$$
Then, the Delta method applied to \eqref{HH_joint_asymptotics} in Lemma~\ref{lemmaH} leads to
\begin{equation*}
%\label{HC_joint_asymptotics}
\sqrt{n}\, \left(\begin{array}{c} H_0\,\big(H(\hat{\bfp}) - H(\bfp)\big)\\ C(\hat{\bfp}) - C(\bfp) \end{array}\right) \xrightarrow{d} N(\mathbf{0}, \Sigma^{(3)})
\text{ with } \Sigma^{(3)}=B_3\,\Sigma^{(1)}\, B_3^\top,
\end{equation*}
where
\begin{eqnarray*}
\sigma^{(3)}_{11} &=& H_0^2\,\sigma^{(1)}_{11},
%\ =\ H_0^2\,\sigma^{(2)}_{11},
\\
\sigma^{(3)}_{12} &=& H_0^2D_0\,\big((D(\bfp) - \tfrac{1}{2}\,H(\bfp))\,\sigma^{(1)}_{11} + H(\bfp)\,\sigma^{(1)}_{12}\big),
%\ =\ H_0^2D_0\,\big(D(\bfp)\,\sigma^{(1)}_{11} + H(\bfp)\,(\sigma^{(1)}_{12} - \tfrac{1}{2}\,\sigma^{(1)}_{11})\big),
%\ =\ H_0^2D_0\,\big(D(\bfp)\,\sigma^{(2)}_{11} + H(\bfp)\,\sigma^{(2)}_{12}\big),
\\
\sigma^{(3)}_{22} &=& H_0^2D_0^2\,\big((D(\bfp) - \tfrac{1}{2}\,H(\bfp))^2\,\sigma^{(1)}_{11} + 2\,H(\bfp)\,(D(\bfp) - \tfrac{1}{2}\,H(\bfp))\,\sigma^{(1)}_{12} + H(\bfp)^2\,\sigma^{(1)}_{22}\big).
%\\
 %&=& H_0^2D_0^2\,\big(D(\bfp)^2\,\sigma^{(2)}_{11} + 2\,H(\bfp)\,D(\bfp)\,\sigma^{(2)}_{12} + H(\bfp)^2\,\sigma^{(2)}_{22}\big).
\end{eqnarray*}
Using the covariances $\sigma^{(2)}_{ij}$ from Theorem~\ref{theorem: asymp distr H D}, we complete the proof of Theorem~\ref{theorem: asymp distr H C}.

\subsection{Proof of Theorem~\ref{theorem: asymp distr H D uniform}}
\label{proof theorem: asymp distr H D uniform}
The second-order Taylor approximation \eqref{PE_Taylor2} adapts to $H\big((\hat{\bfp}+\bfu)/2\big)$ via
$$
H\big((\hat{\bfp}+\bfu)/2\big)\ \approx\ H\big((\bfp+\bfu)/2\big)
 - \sum_{k=1}^{m!} \frac{1+\log{\frac{p_{k}+u_k}{2}}}{2}\,\big(\hat{p}_k-p_{k}\big)
 - \sum_{k=1}^{m!} \frac{\big(\hat{p}_k-p_{k}\big)^2}{4\,(p_{k}+u_k)},
$$
and further reduces for $\bfp=\bfu$ to
$$
\textstyle
H\big((\hat{\bfp}+\bfu)/2\big)\ \approx\ \log{m!}
 - \frac{m!}{8}\,\sum_{k=1}^{m!} \big(\hat{p}_k-p_{k}\big)^2.
$$
Hence, the JS-divergence has the second-order Taylor approximation
\begin{equation}
\label{JS_Taylor2}
D\big(\hat{\bfp}\big)\ \approx\ D(\bfp)
% - \sum_{k=1}^{m!} \frac{\log{\frac{p_{k}+u_k}{2}}-\log{p_{k}}}{2}\,\big(\hat{p}_k-p_{k}\big)
 - \sum_{k=1}^{m!} \frac{\log{\frac{p_{k}+u_k}{2 p_{k}}}}{2}\,\big(\hat{p}_k-p_{k}\big)
 %- \sum_{k=1}^{m!} \Big(\frac{1}{4\,(p_{k}+u_k)} - \frac{1}{4\,p_{k}}\Big)\, \big(\hat{p}_k-p_{k}\big)^2,
 + \sum_{k=1}^{m!} \frac{u_k}{4\,p_{k}(p_{k}+u_k)}\, \big(\hat{p}_k-p_{k}\big)^2,
\end{equation}
which reduces for $\bfp=\bfu$ to
\begin{equation}
\label{JS_Taylor2_uniform}
\textstyle
D\big(\hat{\bfp}\big)\ \approx\ 
 %- (\frac{m!}{8}-\frac{m!}{4})\,\sum_{k=1}^{m!} \big(\hat{p}_k-p_{k}\big)^2
%\ =\ 
 \frac{m!}{8}\,\sum_{k=1}^{m!} \big(\hat{p}_k-p_{k}\big)^2.
\end{equation}
Thus, together with \eqref{PE_Taylor2_uniform} and \eqref{H_QFasymp}, the assertion of Theorem~\ref{theorem: asymp distr H D uniform} follows.

\subsection{Proof of Theorem~\ref{theorem: asymp distr H C uniform}}
\label{proof theorem: asymp distr H C uniform}
We start by deriving a second-order Taylor approximation for~$C(\hat{\bfp})$ from \eqref{eqComplexity}. Using the product rule as well as the above partial derivatives for~$H(\bfp)$ and~$D(\bfp)$, recall \eqref{PE_Taylor2} and \eqref{JS_Taylor2}, we obtain
\begin{eqnarray*}
\tfrac{\partial}{\partial p_{k}}\,C(\bfp)
&=&
H_0 D_0\,\Big(\big(\tfrac{\partial}{\partial p_{k}}\,H(\bfp)\big)\, D(\bfp) + H(\bfp)\, \big(\tfrac{\partial}{\partial p_{k}}\,D(\bfp)\big)\Big)
%\\
%&=&
%H_0 D_0\,\Big(-\big(1+\log{p_{k}}\big)\, D(\bfp) - H(\bfp)\, \frac{\log{\frac{p_{k}+u_k}{2 p_{k}}}}{2}\Big)
\\
&=&
-H_0 D_0\,\Big(\big(1+\log{p_{k}}\big)\, D(\bfp) + \tfrac{1}{2}\,\log{\tfrac{p_{k}+u_k}{2 p_{k}}}\, H(\bfp)\Big),
\end{eqnarray*}
as well as
\begin{eqnarray*}
\tfrac{\partial^2}{\partial p_{k}^2}\,C(\bfp)
&=&
-H_0 D_0\,\Big(\tfrac{1}{p_{k}}\, D(\bfp)
 + \big(1+\log{p_{k}}\big)\, \big(\tfrac{\partial}{\partial p_{k}}\,D(\bfp)\big)
\\
&&\qquad
 %+ \tfrac{1}{2}\,\frac{\tfrac{-u_k}{2 p_{k}^2}}{\tfrac{p_{k}+u_k}{2 p_{k}}}\, H(\bfp)
 + \tfrac{1}{2}\,\tfrac{-u_k}{p_{k}(p_{k}+u_k)}\, H(\bfp)
 + \tfrac{1}{2}\,\log{\tfrac{p_{k}+u_k}{2 p_{k}}}\, \big(\tfrac{\partial}{\partial p_{k}}\,H(\bfp)\big)\Big)
%\\
%&=&
%-H_0 D_0\,\Big(\tfrac{1}{p_{k}}\, D(\bfp)
 %- \tfrac{1}{2}\,\big(1+\log{p_{k}}\big)\, \log{\tfrac{p_{k}+u_k}{2 p_{k}}}
%\\
%&&\qquad
 %- \tfrac{1}{2}\,\tfrac{u_k}{p_{k}(p_{k}+u_k)}\, H(\bfp)
 %- \tfrac{1}{2}\,\log{\tfrac{p_{k}+u_k}{2 p_{k}}}\, \big(1+\log{p_{k}}\big)\Big)
\\
&=&
-H_0 D_0\,\Big(\tfrac{1}{p_{k}}\, D(\bfp)
 - \big(1+\log{p_{k}}\big)\, \log{\tfrac{p_{k}+u_k}{2 p_{k}}}
 - \tfrac{1}{2}\,\tfrac{u_k}{p_{k}(p_{k}+u_k)}\, H(\bfp)\Big)
\end{eqnarray*}
and, for $l\not=k$,
\begin{eqnarray*}
\tfrac{\partial^2}{\partial p_{k}\,\partial p_{l}}\,C(\bfp)
&=&
-H_0 D_0\,\Big(\big(1+\log{p_{k}}\big)\, \big(\tfrac{\partial}{\partial p_{l}}\,D(\bfp)\big)
 + \tfrac{1}{2}\,\log{\tfrac{p_{k}+u_k}{2 p_{k}}}\, \big(\tfrac{\partial}{\partial p_{l}}\,H(\bfp)\big)\Big)
\\
&=&
H_0 D_0\,\Big(\big(1+\log{p_{k}}\big)\, \tfrac{1}{2}\,\log{\tfrac{p_{l}+u_l}{2 p_{l}}}
 + \tfrac{1}{2}\,\log{\tfrac{p_{k}+u_k}{2 p_{k}}}\, \big(1+\log{p_{l}}\big)\Big).
\end{eqnarray*}
If $\bfp=\bfu$, then $H(\bfp)=H_0^{-1}=\log{m!}$, $C(\bfp)=D(\bfp)=0$, and $\log{\tfrac{p_{k}+u_k}{2 p_{k}}}=\log(1)=0$, so the second-order Taylor approximation
\begin{equation}
\label{C_Taylor2_uniform}
\textstyle
C\big(\hat{\bfp}\big)\ \approx\ 
 \frac{m!}{8}\, D_0\,\sum_{k=1}^{m!} 
 \big(\hat{p}_k-p_{k}\big)^2
\end{equation}
follows, and the proof of Theorem~\ref{theorem: asymp distr H C uniform} is completed.

%%%%%%%%%%%%%%%%
\section{MA(q) Process for OPs of Length 3}
\label{section: MA(q) m=3}

In what follows, we establish some results on the covariance matrix in \eqref{eq: long-run covariance matrix} obtained from an MA($q$) process with equal weights for OPs of length $m=3$. 
By \citet[Proposition~7]{bandt07}, the monotone OPs of length 3 have probability $\tfrac{1}{4}$ each, while all the other OPs of length 3 have probability $\tfrac{1}{8}$ each. Therefore,
\begin{equation*}
    \bigl(p_i (\delta_{ij} - p_j)\bigr)_{1\leq i,j \leq m!} = \scalemath{0.9}{\begin{pmatrix}
        \tfrac{3}{16} & \tfrac{1}{32} & \tfrac{1}{32} & \tfrac{1}{32} & \tfrac{1}{32} & \tfrac{1}{16} \\
        \tfrac{1}{32} & \tfrac{7}{64} & \tfrac{1}{64} & \tfrac{1}{64} & \tfrac{1}{64} & \tfrac{1}{32} \\
        \tfrac{1}{32} & \tfrac{1}{64} & \tfrac{7}{64} & \tfrac{1}{64} & \tfrac{1}{64} & \tfrac{1}{32} \\
        \tfrac{1}{32} & \tfrac{1}{64} & \tfrac{1}{64} & \tfrac{7}{64} & \tfrac{1}{64} & \tfrac{1}{32} \\
        \tfrac{1}{32} & \tfrac{1}{64} & \tfrac{1}{64} & \tfrac{1}{64} & \tfrac{7}{64} & \tfrac{1}{32} \\
        \tfrac{1}{16} & \tfrac{1}{32} & \tfrac{1}{32} & \tfrac{1}{32} & \tfrac{1}{32} & \tfrac{3}{16}
    \end{pmatrix}}
    = \frac{1}{64} \scalemath{0.9}{\begin{pmatrix}
        12 & 2 & 2 & 2 & 2 & 4 \\
        2 & 7 & 1 & 1 & 1 & 2 \\
        2 & 1 & 7 & 1 & 1 & 2 \\
        2 & 1 & 1 & 7 & 1 & 2 \\
        2 & 1 & 1 & 1 & 7 & 2 \\
        4 & 2 & 2 & 2 & 2 & 12
    \end{pmatrix}}.
\end{equation*}
Determining the remaining terms is more involved, which results from the fact that we have to consider more than just the increasing and decreasing OPs now. We will make use of the idea of \citet{bandt07} to use a sequence $\delta_1 \delta_2 \delta_3$ of $m=3$ symbols $<$ or $>$. Every OP can be uniquely determined by the relations between every pair of components; for $m=3$ this is $X_1 \delta_1 X_2$, $X_2 \delta_2 X_3$ and $X_1 \delta_3 X_3$. All combinations are enlisted in Table \ref{tab: delta table} in lexicographic order.

\begin{table}[t]
    \label{tab: delta table}
    \caption{Sequences $\delta_1 \delta_2 \delta_3$ for all OPs of length $m=3$.}
    \centering
    \begin{tabular}{cccccc}
    \toprule
        $\pi_1$ & $\pi_2$ & $\pi_3$ & $\pi_4$ & $\pi_5$ & $\pi_6$ \\ 
        \midrule
        (1, 2, 3) & (1, 3, 2) & (2, 1, 3) & (2, 3, 1) & (3, 1, 2) & (3, 2, 1) \\
        $<<<$ & $<><$ & $><<$ & $<>>$ & $><>$ & $>>>$ \\
        \bottomrule
    \end{tabular}
\end{table}

Thus, considering the event that $(X_0, X_1, X_2)$ has OP $\pi_i$ with corresponding sequence $\delta=\delta_1 \delta_2 \delta_3$, it can be rewritten by $X_0 \delta_1 X_1$, $X_1 \delta_2 X_2$ and $X_0 \delta_3 X_2$. Incorporating the definition of $X_t$ in terms of $(\varepsilon_t)_{t\in\bbz}$ leads to 
\begin{equation}
    \label{eq: maq m3 step1}
    \varepsilon_{-q}\delta_1 \varepsilon_1,\, \varepsilon_{1-q} \delta_2 \varepsilon_2 \,\text{ and } \, \varepsilon_{-q} + \varepsilon_{1-q} \delta_3 \varepsilon_1 + \varepsilon_2.
\end{equation}
Assuming $(X_k, X_{k+1}, X_{k+2})$ has OP $\pi_j$ with corresponding sequence $\delta^\prime=\delta^\prime_1 \delta^\prime_2 \delta^\prime_3$ in an analogous way, we obtain
\begin{equation}
    \label{eq: maq m3 step2}
    \varepsilon_{k-q}\delta^\prime_1 \varepsilon_{k+1},\, \varepsilon_{k+1-q} \delta^\prime_2 \varepsilon_{k+2} \,\text{ and } \, \varepsilon_{k-q} + \varepsilon_{k+1-q} \delta^\prime_3 \varepsilon_{k+1} + \varepsilon_{k+2}.
\end{equation}
Note that for $2 \leq k \leq q-1$, there are no intercorrelations between \eqref{eq: maq m3 step1} and \eqref{eq: maq m3 step2}. Hence,
\begin{align*}
    \bbp(\Pi(\bfX_0)=\pi_i, \Pi(\bfX_k)=\pi_j) 
    &= \bbp(\varepsilon_{-q}\delta_1 \varepsilon_1,\ \varepsilon_{1-q} \delta_2 \varepsilon_2, \varepsilon_{-q} + \varepsilon_{1-q} \delta_3 \varepsilon_1 + \varepsilon_2, \\
    &\hspace{10mm}\varepsilon_{k-q}\delta^\prime_1 \varepsilon_{k+1},\ \varepsilon_{k+1-q} \delta^\prime_2 \varepsilon_{k+2}, \varepsilon_{k-q} + \varepsilon_{k+1-q} \delta^\prime_3 \varepsilon_{k+1} + \varepsilon_{k+2}) \\
    &= \bbp(\varepsilon_{-q}\delta_1 \varepsilon_1,\ \varepsilon_{1-q} \delta_2 \varepsilon_2, \varepsilon_{-q} + \varepsilon_{1-q} \delta_3 \varepsilon_1 + \varepsilon_2) \\
    &\hspace{10mm}\cdot\bbp(\varepsilon_{k-q}\delta^\prime_1 \varepsilon_{k+1},\ \varepsilon_{k+1-q} \delta^\prime_2 \varepsilon_{k+2}, \varepsilon_{k-q} + \varepsilon_{k+1-q} \delta^\prime_3 \varepsilon_{k+1} + \varepsilon_{k+2}) \\
    &= \bbp(X_0 \delta_1 X_1, X_1 \delta_2 X_2, X_0 \delta_3 X_2) \cdot \bbp(X_0 \delta^\prime_1 X_1, X_1 \delta^\prime_2 X_2, X_0 \delta^\prime_3 X_2) \\
    &= p_i \cdot p_j
\end{align*}
for $2 \leq k \leq q-1$.
Since this does not depend on the choice of $\delta$ and $\delta^\prime$, the same holds true with reversed roles of $\pi_i$ and $\pi_j$, so the summands for $2 \leq k \leq q-1$ in \eqref{eq: long-run covariance matrix} vanish. We continue by considering the probabilities for $k \in \{1, q, q+1, q+2\}$ separately. Note that for $q=1$, the cases $k=1$ and $k=q$ coincide. 

For $k=1$, we have to consider the cases $q=1$ and $q \geq 2$ separately:
\begin{align}
    &\bbp(\Pi(\bfX_0)=\pi_i, \Pi(\bfX_1)=\pi_j) \nonumber \\
    &\hspace{10mm} = \bbp(\varepsilon_{-q}\delta_1 \varepsilon_1,\ \varepsilon_{1-q} \delta_2 \varepsilon_2,\ \varepsilon_{-q} + \varepsilon_{1-q} \delta_3 \varepsilon_1 + \varepsilon_2,\ \varepsilon_{1-q}\delta^\prime_1 \varepsilon_{2},\ \varepsilon_{2-q} \delta^\prime_2 \varepsilon_{3},\ \varepsilon_{1-q} + \varepsilon_{2-q} \delta^\prime_3 \varepsilon_{2} + \varepsilon_{3}) \nonumber \\
    %&\hspace{20mm} = \indicator{\delta_2=\delta_1^\prime} \cdot \bbp(\varepsilon_{-q}\delta_1 \varepsilon_1,\ \varepsilon_{1-q} \delta_2 \varepsilon_2,\ \varepsilon_{2-q} \delta^\prime_2 \varepsilon_{3},\ \varepsilon_{-q} + \varepsilon_{1-q} \delta_3 \varepsilon_1 + \varepsilon_2,\ \varepsilon_{1-q} + \varepsilon_{2-q} \delta^\prime_3 \varepsilon_{2} + \varepsilon_{3}) \nonumber \\
    &\hspace{10mm} = \indicator{\delta_2=\delta_1^\prime} \cdot \begin{cases}
        \bbp(\varepsilon_{-1}\delta_1 \varepsilon_1 \delta^\prime_2 \varepsilon_{3},\ \varepsilon_0 \delta_2 \varepsilon_2,\ \varepsilon_{-1} + \varepsilon_0 \delta_3 \varepsilon_1 + \varepsilon_2,\ \varepsilon_0 + \varepsilon_1 \delta^\prime_3 \varepsilon_{2} + \varepsilon_{3}), &q=1, \\
        \begin{aligned}[b]
            \bbp(\varepsilon_{-q}\delta_1 \varepsilon_1,\ &\varepsilon_{1-q} \delta_2 \varepsilon_2,\ \varepsilon_{2-q} \delta^\prime_2 \varepsilon_{3},\\ &\varepsilon_{-q} + \varepsilon_{1-q} \delta_3 \varepsilon_1 + \varepsilon_2,\ \varepsilon_{1-q} + \varepsilon_{2-q} \delta^\prime_3 \varepsilon_{2} + \varepsilon_{3}),
        \end{aligned} &q \geq 2,
    \end{cases} \nonumber \\
    &\hspace{10mm} = \indicator{\delta_2=\delta_1^\prime} \cdot \begin{cases}
        \bbp(\varepsilon_1\delta_1 \varepsilon_3 \delta^\prime_2 \varepsilon_{5},\ \varepsilon_2 \delta_2 \varepsilon_4,\ \varepsilon_1 + \varepsilon_2 \delta_3 \varepsilon_3 + \varepsilon_4,\ \varepsilon_2 + \varepsilon_3 \delta^\prime_3 \varepsilon_4 + \varepsilon_5), &q=1, \\
        \bbp(\varepsilon_1\delta_1 \varepsilon_4,\ \varepsilon_2 \delta_2 \varepsilon_5,\ \varepsilon_3 \delta^\prime_2 \varepsilon_6,\ \varepsilon_1 + \varepsilon_2 \delta_3 \varepsilon_4 + \varepsilon_5,\ \varepsilon_2 + \varepsilon_3 \delta^\prime_3 \varepsilon_5 + \varepsilon_6), &q \geq 2. \label{eq: maq m3 1-term}
    \end{cases}
\end{align}

The case $k=q$ and $q \geq 2$ gives
\begin{align}
    \bbp(\Pi(\bfX_0)=\pi_i, \Pi(\bfX_q)=\pi_j) 
    &= \bbp(\varepsilon_{-q}\delta_1 \varepsilon_1,\ \varepsilon_{1-q} \delta_2 \varepsilon_2, \varepsilon_{-q} + \varepsilon_{1-q} \delta_3 \varepsilon_1 + \varepsilon_2, \nonumber \\
    &\hspace{20mm}\varepsilon_0\delta^\prime_1 \varepsilon_{q+1},\ \varepsilon_1 \delta^\prime_2 \varepsilon_{q+2}, \varepsilon_0 + \varepsilon_1 \delta^\prime_3 \varepsilon_{q+1} + \varepsilon_{q+2}) \nonumber \\
    &= \bbp(\varepsilon_1\delta_1\varepsilon_4\delta_2^\prime\varepsilon_7,\ \varepsilon_2\delta_2\varepsilon_5,\ \varepsilon_3\delta_1^\prime\varepsilon_6, \nonumber \\
    &\hspace{20mm} \varepsilon_1+\varepsilon_2\delta_3\varepsilon_4+\varepsilon_5,\ \varepsilon_3+\varepsilon_4\delta_3^\prime\varepsilon_6+\varepsilon_7). \label{eq: maq m3 q-term}
\end{align}
%Let us emphasize that this term is disregarded if $q=1$ and hence, it always holds $2 < q+1$.
In the same manner, for $k=q+1$ and $k=q+2$, we have
\begin{align}
    \bbp(\Pi(\bfX_0)=\pi_i, \Pi(\bfX_{q+1})=\pi_j) 
    &= \bbp(\varepsilon_{-q}\delta_1 \varepsilon_1,\ \varepsilon_{1-q} \delta_2 \varepsilon_2,\ \varepsilon_{-q} + \varepsilon_{1-q} \delta_3 \varepsilon_1 + \varepsilon_2, \nonumber \\
    &\hspace{20mm}\varepsilon_1\delta^\prime_1 \varepsilon_{q+2},\ \varepsilon_{2} \delta^\prime_2 \varepsilon_{q+3},\ \varepsilon_{1} + \varepsilon_{2} \delta^\prime_3 \varepsilon_{q+2} + \varepsilon_{q+3}) \nonumber \\
    &= \bbp(\varepsilon_1\delta_1\varepsilon_3\delta_1^\prime\varepsilon_5,\ \varepsilon_2\delta_2\varepsilon_4\delta_2^\prime\varepsilon_6,\ \varepsilon_1+\varepsilon_2 \delta_3 \varepsilon_3+\varepsilon_4 \delta_3^\prime \varepsilon_5 + \varepsilon_6) \label{eq: maq m3 (q+1)-term}
\end{align}
and
\begin{align}
    \bbp(\Pi(\bfX_0)=\pi_i, \Pi(\bfX_{q+2})=\pi_j) 
    &= \bbp(\varepsilon_{-q}\delta_1 \varepsilon_1,\ \varepsilon_{1-q} \delta_2 \varepsilon_2,\ \varepsilon_{-q} + \varepsilon_{1-q} \delta_3 \varepsilon_1 + \varepsilon_2, \nonumber \\
    &\hspace{20mm}\varepsilon_2\delta^\prime_1 \varepsilon_{q+3},\ \varepsilon_{3} \delta^\prime_2 \varepsilon_{q+4},\ \varepsilon_{2} + \varepsilon_{3} \delta^\prime_3 \varepsilon_{q+3} + \varepsilon_{q+4}) \nonumber \\
    &= \bbp(\varepsilon_1\delta_1\varepsilon_3,\ \varepsilon_2\delta_2\varepsilon_4\delta_1^\prime\varepsilon_6,\ \varepsilon_5\delta_2^\prime\varepsilon_7, \nonumber \\
    &\hspace{20mm} \varepsilon_1 + \varepsilon_2 \delta_3 \varepsilon_3+\varepsilon_4,\ \varepsilon_4+\varepsilon_5\delta_3^\prime\varepsilon_6+\varepsilon_7), \label{eq: maq m3 (q+2)-term}
\end{align}
respectively. Thus, the probabilities \eqref{eq: maq m3 q-term}--\eqref{eq: maq m3 (q+2)-term} are independent of $q$.
It is clear that reversed roles of $\delta$ and $\delta^\prime$ lead to similar results. 

The above considerations reveal that the covariance matrix from \eqref{eq: long-run covariance matrix} remains unaltered for $q\geq2$. However, in contrast to the case $m=2$, the resulting covariance matrix seems to depend on the distribution of the innovations $(\varepsilon_t)_{t\in\bbz}$.

\end{document}